\documentclass[11pt,reqno,twoside]{amsart}

\usepackage{amsmath, amssymb, amscd, amsthm, amscd}

\theoremstyle{plain}
\newtheorem{thm}{Theorem}[section]
\newtheorem{cor}[thm]{Corollary}

\newtheorem{lem}[thm]{Lemma}
\newtheorem{prop}[thm]{Proposition}

\theoremstyle{remark}
\newtheorem{rem}[thm]{Remark}
\newcounter{sspar}[subsection]
\renewcommand\thesspar{(\thesubsection.\arabic{sspar})}
    {\par\ \newline
     \vskip-\baselineskip\vskip.1truecm
     \noindent\refstepcounter{sspar}
     \noindent\textbf{\thesspar} \ignorespaces}
    {\vskip-\baselineskip
    \ignorespaces}
    {\refstepcounter{sspar}
     \textup{\textbf{\thesspar}} \ignorespaces}
    {\vskip-\baselineskip
    \ignorespaces}
\newcommand{\R}{\mathbb R}

\newcommand{\C}{\mathbb C}
\newcommand{\Z}{\mathbb Z}

\newcommand{\al}{\alpha}
\newcommand{\be}{\beta}
\newcommand{\ga}{\gamma}
\newcommand{\Ga}{\Gamma}
\newcommand{\de}{\delta}
\newcommand{\De}{\Delta}
\newcommand{\eps}{\varepsilon}
\newcommand{\si}{\sigma}

\newcommand{\Te}{\Theta}
\newcommand{\la}{\lambda}
\newcommand{\La}{\Lambda}

\newcommand{\inprod}[2]{\langle #1,#2 \rangle}
\newcommand{\F}[5]{\,_{#1}F_{#2} \left( \genfrac{.}{.}{0pt}{}{#3}{#4}
\ ;#5 \right)}
\newcommand{\hf}{\frac{1}{2}}
\newcommand{\Res}[1]{\underset{#1}{\mathrm{Res}}}
\numberwithin{equation}{section}

\begin{document}
\date{\today}
\title{The Wilson function transform}
\author{Wolter Groenevelt}
\address{Technische Universiteit Delft, EWI-TWA \\
Postbus 5031, 2600 GA Delft, The Netherlands}
\email{W.G.M.Groenevelt@ewi.tudelft.nl}
\begin{abstract}
Two unitary integral transforms with a very-well poised $_7F_6$-function as a kernel are given. For both integral transforms the inverse is the same as the original transform after an involution on the parameters. The $_7F_6$-function involved can be considered as a non-polynomial extension of the Wilson polynomial, and is therefore called a Wilson function. The two integral transforms are called a Wilson function transform of type I and type II. Furthermore, a few explicit transformations of hypergeometric functions are calculated, and it is shown that the Wilson function transform of type I maps a basis of orthogonal polynomials onto a similar basis of polynomials.
\end{abstract}
\maketitle

\section{Introduction}
The Jacobi polynomials are orthogonal polynomials that that are solutions to the hypergeometric second order differential equation. The Jacobi polynomials have an explicit expression as $_2F_1$-hypergeometric series, see e.g.~\cite{AAR}. Wilson \cite{Wil} gave a generalization of the Jacobi polynomials as orthogonal $_4F_3$-polynomials. These polynomials are nowadays called the Wilson polynomials. The Wilson polynomials no longer satisfy a second order differential equation, but a second order difference equation. Another way to generalize the Jacobi polynomials is to consider non-polynomial solutions of the hypergeometric differential equation. Spectral analysis of the hypergeometric differential equation leads to a unitary integral transform, called the Jacobi function transform, see e.g.~\cite{Koo}. The kernel in this integral transform is given by a non-polynomial $_2F_1$-functions called the Jacobi function. Koornwinder mentions in \cite[\S9]{Koo2} that there probably also exist functions, naturally called Wilson functions, that generalize the Wilson polynomials as well as the Jacobi functions. 

In this paper we consider non-polynomial eigenfunctions to the Wilson second order difference operator. To study the associated Wilson polynomials, Ismail et al.~\cite{ILVW}, and Masson \cite{Mas}, showed that general solutions of the difference equation are given by very-well-poised $_7F_6$-functions. The analytic part of a certain non-polynomial solution we call a Wilson function. We show that the Wilson functions are the kernel in two unitary integral transforms, which we call the Wilson function transform of type I and type II. The Wilson function satisfies the property, called the duality property, that after an involution on the parameters, the geometric and the spectral parameter can be interchanged. The parameters obtained from the involution are called the dual parameters. The inverse of the Wilson function transform of type I, respectively type II, is again a Wilson function transform of type I, respectively type II, with dual parameters. The transforms can be made completely self-dual by choosing a fixed point of the involution.

The Wilson function transform of type I is a generalized Fourier transform with respect to the same measure as the Wilson polynomials. This measure consists of an absolutely continuous part supported on $[0,\infty)$, and a finite number of discrete mass points. The Wilson function transform of type II is a generalized Fourier transform with respect to a one-parameter family of measures which consist of an absolutely continuous part supported on $[0,\infty)$, and an infinite number of discrete mass points. The extra parameter labels the different sets of discrete mass points.

The Wilson functions can be considered as a formal limit case of the Askey-Wilson functions for $q \rightarrow 1$. The Askey-Wilson functions are eigenfunction to the Askey-Wilson second order difference operator for $0 < q  <1$. The Askey Wilson functions are the kernel in the Askey-Wilson function transform, which is found by Koelink and Stokman \cite{KSt}. They show in \cite{KSt2} that the Askey-Wilson function has an interpretation as spherical function for the quantum $SU(1,1)$ group. From this point of the view the Askey-Wilson function is a $q$-analogue of the Jacobi function, since the Jacobi function has similar interpretation for the Lie group $SU(1,1)$, see \cite{Koo}. So the Wilson functions in this paper give a new limit case of the Askey-Wilson functions. We do not know if the Wilson function also has an interpretation as a spherical function. 

The Wilson function can also be considered as a formal limit case of Ruijsenaars' $R$-function \cite{Ru1}, \cite{Ru2}, \cite{Ru3}, which is an eigenfunction to the Askey-Wilson second order difference operator for $|q|=1$. The $R$-function is given by a Barnes-type integral which is considered as a generalization of the Barnes integral representation for the $_2F_1$-series. Using the Barnes-type integral representation for the $_7F_6$-series \cite[(4.7.1.3)]{Sl1}, the $R$-function can also be considered as a generalization of the Wilson function.  

In a future paper we will show that the Wilson functions, and also the Wilson polynomials, have an interpretation as Racah coefficients for tensor products of positive discrete series, negative discrete series, and principal unitary series representations of the Lie algebra $\mathfrak{su}(1,1)$. Both Wilson function transforms in this paper have an interpretation in the context of Racah coefficients.

The organization of this paper is as follows. In section \ref{sec:Wpol} we give some well-known properties of the Wilson polynomials. The Wilson polynomials are eigenfunctions of a second order difference operator $\La$, and we show that the Wilson polynomials are also eigenfunction of the same difference operator with dual parameters.

In section \ref{sec:Wilson f} we consider a certain type of non-polynomial eigenfunction of a difference operator $L$ which is closely related to $\La$. These eigenfunctions, the Wilson functions, are also eigenfunctions of the difference operator $L$ with dual parameters. 

In section \ref{sec:I} we define a Hilbert space $\mathcal M$, and, using the asymptotic behaviour of the Wilson function, we show that a truncated inner product of two Wilson functions approximates a reproducing kernel. This leads to a unitary integral transform, which we call the Wilson function transform of type I. 

In section \ref{sec II} we define a different Hilbert space $\mathcal H$. Again using asymptotic behaviour of the Wilson function, we show that a truncated inner product of two Wilson functions approximates a reproducing kernel. This leads to the Wilson function transform of type II.

In section \ref{sec:expl} the Wilson function transforms of a Jacobi function, and of a Wilson polynomial, are calculated explicitly, using an integral representation of the Wilson function. Also we show that the Wilson function transform of type I maps an orthogonal system of Wilson polynomials onto the same orthogonal system with dual parameters. \\

\emph{Notations.}
We use the standard notation for the hypergeometric series, i.e.
\[
\F{p}{q}{a_1, \ldots, a_{p}} {b_1, \ldots, b_q}{z} = \sum_{n=0}^\infty \frac{ (a_1)_n \ldots (a_{p})_n }{ (b_1)_n \ldots (b_q)_n } \frac{z^n}{n!},
\]
where $(a)_n$ denotes the Pochhammer symbol, defined by
\[
(a)_n = \frac{\Ga(a+n)}{\Ga(a)}=a(a+1)(a+2) \ldots (a+n-1), \qquad n \in \Z_{\geq 0}.
\]
A hypergeometric series is called very-well poised if $p=q+1$, $a_1+1 = b_1+a_2 = \ldots = b_q+a_{q+1}$ and $b_1 = a_2/2$. For a very-well poised $_7F_6$-series of argument $1$ we use Bailey's $W$ notation, see \cite{Bai}, i.e.
\[
W(a;b,c,d,e,f) = \F{7}{6}{a,\ 1+\hf a,\ b,\ c,\ d,\ e,\ f\ }{\hf a, 1+a-b, 1+a-c, 1+a-d, 1+a-e, 1+a-f}{1}.
\]
If the series does not terminate, the condition for convergence is $\Re(2a+2-b-c-d-e-f)>0$.
Also, since the product $\Ga(a+b) \Ga(a-b)$ frequently occurs in this paper, we use for this product the shorthand notation $\Ga(a \pm b)$.

\section{Wilson polynomials} \label{sec:Wpol}
In this section we recall some well-known properties of the Wilson polynomials.  The Wilson polynomial satisfies a second order difference equation in its degree (the three-term recurrence relation), and also a second order difference equation in its argument. The goal of this section is to point out that the two difference equations are the same after a change of the parameters.\\

The Wilson polynomials $R_n(x)$, see \cite{Wil}, \cite[\S6.10]{AAR}, \cite{KS}, are polynomials in $x^2$ of degree $n$. They can be defined by the initial values $R_{-1}(x)=0$, $R_0(x)=1$, and the recurrence relation
\begin{equation} \label{eq:rec Wpol}
-(a^2+x^2) R_n(x) = C_n \big[ R_{n+1}(x)-R_{n}(x)\big] + D_n \big[R_{n-1}(x) - R_n(x)\big],
\end{equation}
where
\[
\begin{split}
C_n &= C_n(a,b,c,d) =\frac{ (n+a+b+c+d-1)(n+a+b)(n+a+c)(n+a+d) }{ (2n+a+b+c+d-1)(2n+a+b+c+d) } , \\
D_n &= D_n(a,b,c,d) =\frac{ n(n+b+c-1)(n+b+d-1)(n+c+d-1) }{ (2n+a+b+c+d-2)(2n+a+b+c+d-1) }.
\end{split}
\]
The explicit expression for the polynomials $R_n$ is given by
\begin{equation} \label{def:Wilson pol}
R_n(x)= R_n(x;a,b,c,d) = \F{4}{3}{-n, n+a+b+c+d-1, a+ix, a-ix}{a+b, a+c, a+d}{1}.
\end{equation}

Let $a,b,c,d \in \C$ be such that non-real parameters appear in conjugate pairs with positive real part, and such that the pairwise sum of any two parameters has positive real part. We define the measure $d\mu(x)=d\mu(x;a,b,c,d)$ by
\[
\int f(x) d\mu(x) =\frac{1}{2\pi} \int_0^\infty f(x) w(x) dx + i\sum_k f(x_k) w_k,
\]
where
\[
w(x) = w(x;a,b,c,d) =  \frac{ \Ga(a\pm ix) \Ga(b\pm ix) \Ga(c\pm ix) \Ga(d\pm ix) }{ \Ga(\pm 2ix) }.
\] 
The points $x_k$ are of the form $x_k=i(e+k)$, where $e$ is any of the parameters $a,b,c,d$ with $e<0$. The sum is over $k \in \Z_{\geq 0}$, such that $e+k<0$. The weights $w_k$ are the residue at $x=x_k$ of $w(x)$. In particular, if all parameters are positive, or occur in pairs of complex conjugates with positive real part, the measure $d\mu$ is absolutely continuous. The Wilson polynomials  are orthogonal with respect to the measure $d\mu$, i.e.
\begin{multline} \label{eq:orth}
\int R_m(x) R_n(x) d\mu(x) =\\ \de_{nm}\frac{a+b+c+d-1}{2n+a+b+c+d-1}\, \frac{n!\,  \Ga(a+b) \Ga(a+c) \Ga(a+d) \Ga(b+c+n) \Ga(b+d+n) \Ga(c+d+n) }{(a+b)_n (a+c)_n (a+d)_n (a+b+c+d-1)_n\,  \Ga(a+b+c+d)}.
\end{multline}

The Wilson polynomials also satisfy a difference equation in $x$, given by
\begin{equation} \label{eq:diff}
n(n+a+b+c+d-1) R_n(x) = A(-x)\big[ R_n(x+i)- R_n(x) \big] + A(x)\big[ R_n(x-i)- R_n(x) \big],
\end{equation}
where
\[
A(x) = \frac{ (a+ix)(b+ix)(c+ix)(d+ix) }{2ix(2ix+1)}.
\]

We define the difference operator $\La$ by
\begin{equation} \label{def:La}
\La = A(x)(T_{-i}-I) + A(-x)(T_{i}-I),
\end{equation}
where $I$ denotes the identity operator and $T$ is the shift operator (i.e. $T_z f(x) =f(x+z)$). From the difference equation \eqref{eq:diff} it follows that the polynomials $R_n(x;a,b,c,d)$ are eigenfunctions of $\La$ for eigenvalue $n(n+a+b+c+d-1)$. 

The recurrence relation can be written in a self-dual way. Given the parameters $a,b,c,d \in \C$, we define dual parameters by
\begin{equation} \label{eq:dual param}
\begin{split}
\tilde a &= \hf(a+b+c+d-1),\quad 
\tilde b=\hf(a+b-c-d+1) , \\
\tilde c&=\hf(a-b+c-d+1), \quad
\tilde d=\hf(a-b-c+d+1).
\end{split}
\end{equation}
It is an easy verification that $(a,b,c,d) \mapsto (\tilde a, \tilde b, \tilde c, \tilde d)$ defines an involution on the parameters. For a function $f=f(a,b,c,d)$, we define $\tilde f$ by $\tilde f = f(\tilde a, \tilde b, \tilde c, \tilde d)$. We use the same notation for other objects, like measures, sets and operators. Denote the Wilson polynomial $R_n(x)$ by $P_\la(x)$, where $\la = i(n+\tilde a)$, i.e.
\[
P_\la(x) = \F{4}{3}{\tilde a + i\la, \tilde a - i\la, a+ix, a-ix}{a+b, a+c, a+d}{1},
\]
then we see that $P_\la(x) = \tilde P_x(\la)$, since $a+e=\tilde a + \tilde e$, for $e=b,c,d$. The recurrence relation \eqref{eq:rec Wpol} and the difference equation \eqref{eq:diff} for the Wilson polynomials can now be written as 
\[
\begin{split}
\tilde \La \tilde P_x(\la) &= -(a^2+x^2) \tilde P_x(\la),\\
\La P_\la(x) &= -(\tilde a^2 + \la^2)P_\la(x).
\end{split}
\]
So we see that, for $\la \in i(\tilde a + \Z_{\geq 0})$, the Wilson polynomial $P_\la(x)$ is an eigenfunction of both $\La$ and $\tilde \La$. 
 
\section{Wilson functions} \label{sec:Wilson f}
In the previous section we observed that, for $\la^2 = -(\tilde a+n)^2$, the Wilson polynomials are solutions to the eigenvalue equation 
\begin{equation} \label{eq:eigenv La}
(\La f)(x) = -(\la^2+\tilde a^2)f(x),
\end{equation}
For more general values of $\la$, solutions to \eqref{eq:eigenv La} can be given in terms of very-well poised $_7F_6$-series. This is shown by Ismail et al.~ \cite{ILVW} and by Masson \cite{Mas}, who investigate the associated Wilson polynomials. Let $\psi_\la(x)$ be the function defined by
\[
\begin{split}
\psi_\la(x)=\psi_\la(x;a,b,c,d) =& \frac{ \Ga(b+c) \Ga(\tilde a+\tilde b+\tilde c+i\la) \Ga(1-\tilde d+i\la) \Ga(1-d \pm ix) }{   \Ga(\tilde b+c+i\la \pm ix) }\\
& \times W(\tilde a+\tilde b+ \tilde c-1 +i\la;a+ix, a-ix, \tilde a+i\la, \tilde b+i\la, \tilde c+i\la ).
\end{split}
\]
The dual parameters $\tilde a, \tilde b, \tilde c, \tilde d$ are still defined by \eqref{eq:dual param}. By \cite[Thm.~2]{ILVW}, or \cite[(2.5)]{Mas}, the function $\psi_\la(x)$ is a solution to the eigenvalue equation \eqref{eq:eigenv La}.

Instead of $\psi_\la(x)$, we study the closely related function 
\begin{equation}\label{def:Wilson function}
\begin{split}
\phi_\la(x)=  \phi_\la(x;a,b,c,d) 
=& \frac{ \Ga(\tilde a+\tilde b+\tilde c+i\la) }{\Ga(a+b) \Ga(a+c) \Ga(1+a-d) \Ga(1-\tilde d -i\la)   \Ga(\tilde b+c+i\la \pm ix) }\\
& \times W(\tilde a+\tilde b+ \tilde c-1 +i\la;a+ix, a-ix, \tilde a+i\la, \tilde b+i\la, \tilde c+i\la ).
\end{split}
\end{equation} 
So $\phi_\la(x) = \psi_\la(x)/ K(x,\la)$, where $K(x,\la)$ is the function given by
\[
K(x,\la) = \Ga(a+b) \Ga(a+c) \Ga(b+c) \Ga(1+a-d) \Ga(1-d \pm ix) \Ga(1-\tilde d \pm i\la).
\]
We call the function $\phi_\la(x)$ a Wilson function. The Wilson function $\phi_\la(x)$ has the advantage that it is symmetric in $a,b,c,1-d$, cf.~ Remark \ref{rem:Phi}(iii), while $\psi_\la(x)$ is not. Also $\phi_\la(x)$ is an analytic function in $(x,\la) \in \C^2$.

The $_7F_6$-series in the definition of $\phi_\la(x)$ converges absolutely for $\Re(1-\tilde d - i\la)>0$. Writing the $_7F_6$-series as a sum of two balanced $_4F_3$-series, we have an expression for $\phi_\la(x)$ which always converges: 
\begin{equation} \label{eq:expansion1}
\begin{split}
\phi_\la(x)=& \frac{ \Ga(1-a-d)}{ \Ga(a+b) \Ga(a+c) \Ga(1-d \pm ix) \Ga(1-\tilde d \pm i\la)} \F{4}{3}{ a+ix, a-ix, \tilde a+i\la, \tilde a-i\la }{ a+b, a+c, a+d}{1} \\
+&\, \frac{ \Ga(a+d-1)  }{\Ga(1+b-d) \Ga(1+c-d) \Ga(a \pm ix) \Ga(\tilde a \pm i\la) } \\
&\times \F{4}{3}{1-d+ix, 1-d-ix, 1-\tilde d+i\la, 1-\tilde d-i\la}{1+b-d, 1+c-d, 2-a-d}{1}.
\end{split}
\end{equation}
This follows from \cite[\S4.4(4)]{Bai} with parameters specified by
\[
a \mapsto \tilde a + \tilde b + \tilde c -1 + i\la, \quad c \mapsto \tilde c +i\la, \quad d \mapsto \tilde b +i\la,\quad e \mapsto \tilde a + i\la, \quad f \mapsto a-ix, \quad g \mapsto a+ix.
\]
From \eqref{eq:expansion1} we see that $\phi_\la(x)$ is an analytic function in $(x,\la) \in \C^2$. Observe that for $\la = \pm i(\tilde a + n)$, the second term in \eqref{eq:expansion1} vanishes because of the factor $\Ga(\tilde a \pm i\la)^{-1}$, and then we see that $\psi_\la(x)= K(x,\la)\phi_\la(x)$ reduces to a Wilson polynomial. So $\psi_\la(x)$ is the analytic continuation of the Wilson polynomial in its degree.

From \eqref{eq:expansion1} and \eqref{eq:dual param} follows the duality property 
\begin{equation} \label{eq:duality}
\phi_\la(x) = \tilde \phi_x(\la).
\end{equation}
This duality property is similar to the duality property for the Askey-Wilson functions in \cite{KSt}. For the $_7F_6$-series in the definition of the Wilson function, \eqref{eq:duality} is implied by Bailey's transformation \cite[\S7.5(1)]{Bai}.

Since $\phi_\la(x) = \psi_\la(x)/K(x,\la)$ and $\psi_\la(x)$ is a solution to eigenvalue equation \eqref{eq:eigenv La}, the Wilson function satisfies the equation
\[
M^{-1}_{K(x,\la)} \circ \La \circ M_{K(x,\la)}\ \phi_\la(x) = -(\la^2 + \tilde a^2) \phi_\la(x).
\]
Here $M_{K(x,\la)}$ denotes multiplication by $K(x,\la)$. From this we obtain the following proposition.
\begin{prop} \label{prop:L}
The Wilson function $\phi_\la(x)$ is a solution to 
\[
(L f)(x) = (\la^2+\tilde a^2)f(x),
\]
where $L$ is the difference operator defined by
\[
L = B(-x) T_i + \big[A(-x)+A(x)\big] I + B(x)T_{-i},
\]
\[
\begin{split}
A(x) &= \frac{ (a+ix)(b+ix)(c+ix)(d+ix) }{2ix(2ix+1)},\\
B(x) &= \frac{ (a+ix)(b+ix)(c+ix)(1-d+ix) }{2ix(2ix+1)}.
\end{split}
\]
\end{prop}

In the next two sections we consider the action of the second order difference operator $L$ on two different Hilbert spaces, which leads to two different integral transforms with the Wilson function as a kernel. The method we use comes down to approximating with the Fourier transform, using asymptotic expansions of the Wilson functions. This method is essentially the same method as used by G\"otze \cite{Go}, and Braaksma and Meulenbeld \cite{BM}, for the Jacobi function transform. A similar method is also used by Koelink and Stokman in \cite{KSt} for the Askey-Wilson function transform, and in \cite{GKR} for the continuous Hahn function transform.

\section{The Wilson function transform: type I} \label{sec:I}
\subsection{The Hilbert space $\mathbf{\mathcal M}$}
Let $V\subset \C^4$ be the set of parameters $(a,b,c,d)$ satisfying the following two conditions:
\begin{itemize}
\item The parameters $a,b,c,1-d$ are real except for pairs of complex conjugates with positive real part.
\item The pairwise sum of $a,b,c,1-d$ is contained in $\C \setminus (-\infty,0]$.
\end{itemize}
A direct verification shows that the assignment $(a,b,c,d) \mapsto (\tilde a, \tilde b, \tilde c, \tilde d)$, see \eqref{eq:dual param}, is an involution on $V$. Throughout this section we assume $(a,b,c,d) \in V$.

Let $M$ be the weight given by
\[
M(x)= M(x;a,b,c,d) = \frac{ \Ga(a \pm ix) \Ga(b \pm ix) \Ga(c \pm ix) \Ga(1-d \pm ix)}{  \Ga(\pm 2ix)}.
\]
The weight $M$ is positive for $x \in \R$. Observe that $M(\cdot;a,b,c,d) = w(\cdot;a,b,c,1-d)$, where $w$ is the weight function for the Wilson polynomials. We assume that the function $M$ has only simple poles. This imposes conditions on the parameters $a,b,c,d$ that can be removed afterwards by continuity in the parameters. For $e \in \C$ define the set $\mathcal D_e$ by
\[
\mathcal D_e = \{ i(e+n) \ | \ n \in \Z_{\geq 0},\ e+n< 0\},
\]
and let $\mathcal D = \mathcal D_a \cup \mathcal D_b \cup \mathcal D_c \cup \mathcal D_{1-d}$. 
We define the measure $dm(\cdot)=dm(\cdot;a,b,c,d)$ by
\[
\int f(x) dm(x) = \frac{1}{2\pi}\int_0^\infty f(x) M(x) dx +i \sum_{x \in \mathcal D}\ f(x) \Res{z=x}\, M(z).
\]
If $x \in \mathcal D_a$, we have explicitly
\[
\begin{split}
i\Res{z=i(a+n)}\, M(z) =&
\frac{ \Ga(a+b) \Ga(a+c) \Ga(1+a-d) \Ga(b-a) \Ga(c-a) \Ga(1-d-a) }{ \Ga(-2a) } \\ & \times \frac{(2a)_n (a+1)_n(a+b)_n (a+c)_n (1+a-d)_n }{n!\,(a)_n (1+a-b)_n (1+a-c)_n (a+d)_n },
\end{split}
\]
and then we see that for $(a,b,c,d) \in V$, the measure $dm$ is positive.
Recall that $\phi_\la(x)$ is symmetric in $a,b,c,1-d$, and by \eqref{eq:expansion1} the Wilson function is obviously even in $x$ and $\la$, therefore $\overline{\phi_\la(x)} = \phi_{\overline{\la}}(\overline{x})$. From this it follows that for $x\in \mathrm{supp}\ dm(\cdot)$ and $\la \in \mathrm{supp}\ d\tilde m(\cdot)$ the Wilson function $\phi_\la(x)$ is real valued.

We define the Hilbert space $\mathcal M=\mathcal M(a,b,c,d)$ to be the Hilbert space consisting of even functions that have finite norm with respect to the inner product $\inprod{\cdot}{\cdot}_{\mathcal M}$ defined by
\[
\inprod{f}{g}_{\mathcal M} = \int f(x) \overline{g(x)} dm(x).
\]

\subsection{The Wronskian}
For $0<N<\infty$, we define a pairing $\langle \cdot,\cdot \rangle_{N}$ by
\[
\inprod{f}{g}_{N}= \frac{1}{2\pi}\int_0^N f(x)g(x) M(x)dx +i\ \sum_{x \in \mathcal D}\ f(x)g(x) \Res{z=x}\, M(z).
\]
If $f$ and $g$ are real valued functions in $\mathcal M$, the limit $N \rightarrow \infty$ gives the inner product $\inprod{f}{g}_{\mathcal M}$. \\
For functions $f,g$ that are analytic in $\C$, we define the Wronskian $[f,g]$ by
\[
[f,g](z)=\frac{1}{2\pi}\int_z^{z+i} \left\{ f(x) g(x-i)- f(x-i)g(x) \right\}
B(x) M(x) dx.
\]
\begin{lem} \label{lem:oddWr}
Let $f,g$ be analytic in $\C$ and even, then $z \mapsto [f,g](z)$ is odd in $z$.
\end{lem}
\begin{proof}
Let $I$ be the function given by
\[
I(x) = \Big(f(x) g(x-i) - f(x-i) g(x) \Big) B(x) M(x),
\]
then $[f,g](z)= \int_z^{z+i} I(x)dx$.
Since $f(x),g(x)$ and $M(x)$ are even functions in $x$, and $B(-x) M(x) = B(x+i) M(x+i)$, we have $I(-x) = - I(x+i)$. Therefore
\[
\int_z^{z+i} I(x) dx = - \int_{-z}^{-z-i} I(-x) dx = \int_{-z}^{-z-i} I(x+i) dx = -\int_{-z}^{-z+i} I(x) dx.
\]
Hence $z \mapsto [f,g](z)$ is an odd function in $z$.
\end{proof}
\begin{prop} \label{prop: Wr}
For $N \gg 0$ and for even analytic functions $f$ and $g$ we have
\[
\langle Lf ,  g \rangle_N - \langle f ,  Lg \rangle_N = [f,g](N).
\]
\end{prop}
\begin{proof}
Recall that we assume that the poles of $M(x)$ are simple. For even functions $f$ and $g$ we have
\[
\langle f ,  g \rangle_N = \frac{1}{4\pi}\int_{\mathcal C_N} f(x) g(x) M(x) dx,
\]
where $\mathcal C_N$ is a contour in the complex plane defined as follows:
\begin{itemize}
\item $\mathcal C_N$ starts at $x=-N$ and ends at $x=N$,
\item $\mathcal C_N$ is invariant under reflection in the origin,
\item $\mathcal C_N$ separates the sequence of poles $i(a+n)$, $n \in \Z_{\geq 0}$ from the sequence $-i(a+m)$, $m \in \Z_{\geq 0}$, and similarly for poles of $M(x)$ corresponding to $b,c,1-d$.  
\end{itemize}
Now we have
\[
\begin{split}
\frac{1}{4\pi}\int_{\mathcal C_N}& (Lf)(x) g(x) M(x) dx - \frac{1}{4\pi} \int_{\mathcal C_N} f(x) (Lg)(x) M(x) dx = \\
&\frac{1}{4\pi}\int_{\mathcal C_N} \Big( f(x+i) g(x)  - f(x) g(x+i) \Big) B(-x) M(x) dx \\+& \frac{1}{4\pi}\int_{\mathcal C_N} \Big( f(x-i) g(x) - f(x) g(x-i) \Big) B(x) M(x) dx.
\end{split}
\]
Using $B(-x+i)M(x-i) = B(x)M(x)$, we can write this as
\[
\frac{1}{4\pi}\left(\int_{\mathcal C_N} -\int_{\mathcal C_N+i} \right) \Big( f(x-i) g(x) - f(x) g(x-i) \Big) B(x) M(x) dx.
\]
Here $\mathcal C_N+i= \{ x \in \C \, |\, x-i \in \mathcal C_N \}$. The integrand has its poles at $x=i(e+1+n)$, $x=-i(e+m)$, for $n,m \in \Z_{\geq 0}$ and $e= a,b,c,1-d$. Now we make a closed contour by connecting $\mathcal C_N$ and $\mathcal C_N+i$ at the end points in a straight line (there are no poles of $M$ on these lines for $N$ large enough), then the integrand $I(x)$ has no poles inside the closed contour. So, by Cauchy's Theorem,
\[
\int_{\mathcal C_N} -\int_{\mathcal C_N+i} = \int_{-N}^{-N+i} - \int_N^{N+i},
\]
and from this we obtain
\[
\langle Lf , g \rangle_N - \langle f , Lg \rangle_N = \hf[f,g](N) - \hf[f,g](-N).
\]
Now the proposition follows from Lemma \ref{lem:oddWr}.
\end{proof}
Since the Wilson functions are eigenfunctions of $L$ for eigenvalue $\tilde a^2+ \la^2$, we find from Proposition \ref{prop: Wr} the following.
\begin{prop} \label{prop:int Wr}
For $\la \neq \overline{\la'}$,
\[
\inprod{\phi_\la}{ \phi_{\la'}}_N = \frac{ [\phi_\la, \phi_{\la'}](N)}{\la^2-\la'^2 }. 
\]
\end{prop}

Next we want to let $N \rightarrow \infty$ in Proposition \ref{prop:int Wr}, so we need the asymptotic behaviour of the Wilson function and of $B(x+iy)M(x+iy)$ for $x \rightarrow \infty$ and $0 \leq y \leq 1$. We have
\begin{equation} \label{eq:asym M}
\begin{split}
M(x+iy)=& 
 16\pi^3 x^{2a+2b+2c-2d-1}e^{-2\pi x -2i\pi y}\left(1+ \mathcal O\Big(\frac{1}{x}\Big) \right),
\\
B(x+iy) = & -\frac{x^2}{4}\left(1+ \mathcal O\Big(\frac{1}{x}\Big) \right).
\end{split}
\end{equation}
The asymptotic behaviour of the weight function $M$ can be obtained from 
\cite[\S4.5]{Olv}
\begin{equation} \label{eq:asym Ga}
\frac{ \Ga (a+z)}{ \Ga(b+z) } = z^{a-b}\left( 1+
\frac{1}{2z}(a-b)(a+b-1) + \mathcal O \big(\frac{1}{z^2}\big) \right), \quad |z|\rightarrow \infty, \quad |\arg(z)|<\pi,
\end{equation}
and from applying Euler's reflection formula for the $\Ga$-function
\[
\Ga(a \pm ix) = \frac{\pi \Ga(a+ ix)}{\Ga(1-a+ix) \sin \pi(a-ix)}.
\]
To determine the asymptotic behaviour of $\phi_\la(x)$, we expand $\phi_\la(x)$ in functions with nice asymptotic behaviour.
\begin{prop} \label{prop:c-expa}
\[
\phi_\la(x) = \tilde c(\la) \Phi_\la(x) + \tilde c(-\la) \Phi_{-\la}(x),
\] 
where
\[
\begin{split}
\Phi_\la(x) &= \frac{ 1}{ \Ga( \tilde b+c+i\la \pm ix) } \F{4}{3}{\tilde a+i\la, \tilde b+i\la, \tilde c+i\la, 1-\tilde d+i\la}{ \tilde b+c+i\la+ix, \tilde b+c+i\la -ix, 1+2i\la}{1},\\
\tilde c(\la) &=  \frac{ \Ga(-2i\la) }{ \Ga(\tilde a-i\la) \Ga(\tilde b-i\la) \Ga(\tilde c-i\la) \Ga(1-\tilde d-i\la)}.
\end{split}
\]
\end{prop}
\begin{proof}
This follows from transforming the $_7F_6$-function in the definition \eqref{def:Wilson function} of $\phi_x(\la)$  by \cite[\S4.4(4)]{Bai} with parameters specified by
\[
a \mapsto \tilde a + \tilde b + \tilde c -1 + i\la, \quad c \mapsto a +ix, \quad d \mapsto a-ix,\quad e \mapsto \tilde a + i\la, \quad f \mapsto \tilde b + i\la, \quad g \mapsto \tilde c + i\la.
\] 
\end{proof}

\begin{rem} \label{rem:Phi}
(i) Note that the dual weight function $\tilde M$ can be expressed in terms of the $\tilde c$-function as
\[
\tilde M(\la) = \frac{1}{\tilde c(\la) \tilde c(-\la)}.
\]

(ii) The functions $\Phi_\la$ and $\Phi_{-\la}$ are in general not solutions to the eigenvalue equation in Proposition \ref{prop:L}. Let us define
\[
\begin{split}
&\Psi_\la(x) =  \\
&\frac{ \Ga( 1-a-ix) 
 \Ga( 1-b-ix)  \Ga( 1-c-ix) \Ga( 1-d-ix) \Ga(2-d-ix+2i\la) }{ \Ga(1+ \tilde a-d+i\la-ix) \Ga(1+\tilde d-d+i\la-ix) \Ga(2-\tilde b-d+i\la-ix) \Ga(2-\tilde c-d+i\la-ix) \Ga(d-ix) } \\
&\times \frac{ \sin \pi (1+ \tilde a-d+ix+i\la) }{ \sin \pi (d-ix) }\,W(1-d-ix+2i\la; 1-d -ix, 1-\tilde a+i\la , 1-\tilde d+i\la, \tilde b+i\la, \tilde c+i\la).
\end{split}
\] 
From \cite[(2.12)]{Mas} it follows that $\Psi_\la$ and $\Psi_{-\la}$ are solutions to the eigenvalue equation \eqref{eq:eigenv La}, so $\Psi_\la(x)/K(x,\la)$ is an eigenfunction of $L$ for eigenvalue $\la^2+ \tilde a^2$. The Wilson function $\phi_\la(x)$ can be expanded in terms of $\Psi_\la(x)$ and $\Psi_{-\la}(x)$ as follows:
\[
\Ga(1-d \pm ix) \phi_\la(x) = \tilde c(\la)  \Psi_\la(x) +\tilde c(-\la) \Psi_{-\la}(x),
\]
where the $\tilde c$-function is the same as in Proposition \ref{prop:c-expa}. This expansion follows from \cite[(4.3.7.8)]{Sl1}, with parameters given by
\begin{gather*}
a \mapsto 1-d+2i\la-ix, \quad b \mapsto 1-d-ix, \quad c \mapsto 1-\tilde a+i\la,\\
d \mapsto 1-\tilde d+i\la, \quad e \mapsto \tilde b+i\la, \quad f \mapsto \tilde c+i\la,
\end{gather*}
and Euler's reflection formula for the $\Ga$-function. From this expansion is it possible to determine the asymptotic behaviour of $\phi_\la(x)$ for $x \rightarrow \infty$. We will use the simpler functions $\Phi_{\pm \la}$ from Proposition \ref{prop:c-expa} to compute the asymptotic behaviour of the Wilson function.

(iii) Writing $c = \hf(\tilde a - \tilde b + \tilde c - \tilde d+1)$ in Proposition \ref{prop:c-expa} we see that the Wilson function $\phi_\la(x;a,b,c,d)$ is symmetric in $\tilde a, \tilde b, \tilde c, 1-\tilde d$, and by the duality property \eqref{eq:duality} then also in $a,b,c,1-d$.
\end{rem}

To determine the asymptotic behaviour of $\Phi_\la(x)$ we use Euler's reflection formula to rewrite the $\Ga$-functions in front of the $_4F_3$-function and we use \eqref{eq:asym Ga}. Then we find for $\Phi_\la$, for $y \in \R$ and $x \rightarrow \infty$,
\begin{equation} \label{eq:asympt}
\Phi_\la(x+iy) = \frac{1}{2i\pi} x^{d-a-b-c-2i\la} e^{\pi x + i\pi y}\left( 1 + \frac{y}{ix}(a+b+c-d +2i\la) +\mathcal O\big(\frac{1}{x^2} \Big) \right).
\end{equation}
From Proposition \ref{prop:c-expa} it follows that the asymptotic behaviour of the Wronskian $[\phi_\la, \phi_{\la'}](N)$ can be computed from the four Wronskians $[\Phi_{\pm \la}, \Phi_{ \pm \la'}](N)$. So we must compute $[\Phi_\la, \Phi_{\la'}](N)$ for $N \rightarrow \infty$.
\begin{lem} \label{lem:asympt Wr}
For $N \rightarrow \infty$,
\[
[\Phi_\la,\Phi_{\la'}](N)= i( \la - \la')N^{-2i(\la+\la')}\left(1+ \mathcal O\Big(\frac{1}{N} \Big) \right).
\]
\end{lem}
\begin{proof}
Let $G(x)$ be the function given by
\[
G(x) = \Phi_\la(x) \Phi_{\la'}(x-i) - \Phi_\la(x-i) \Phi_{\la'}(x).
\]
From the asymptotic behaviour of $\Phi_\la$ given above we find, for $0 \leq y \leq 1$ and $x \rightarrow \infty$,
\[
G(x+iy) =-\frac{( \la - \la')}{2\pi^2} x^{2d-2a-2b-2c-1-2i(\la+\la')}e^{2\pi x +2i\pi y}\left( 1+ \mathcal O\Big( \frac{1}{x} \Big) \right). 
\]
Using the asymptotic behaviour of $B(x+iy)M(x+iy)$ for $x \rightarrow \infty$ gives 
\[
G(x+iy)B(x+iy)M(x+iy) =2 \pi( \la -\la') x^{-2i(\la+\la)} \left( 1+ \mathcal O \Big( \frac 1x \Big) \right). 
\]
Note that the main term in de asymptotic expansion is independent of $y$. Next we write
\[
[\Phi_\la,\Phi_{\la'}](N) = \frac{1}{2\pi}\int_N^{N+i} G(x) B(x) M(x) dx =  -\frac{1}{2\pi i} \int_0^1 G(N+iy) B(N+iy) M(N+iy) dy,
\] 
and we apply dominated convergence to obtain the result.
\end{proof}

\subsection{Continuous spectrum}

In this subsection we assume $\la, \la' \in \R$, and since $\phi_\la$ is even in $\la$, we may assume $\la, \la' \geq 0$.

\begin{prop} \label{prop1}
Let $f$ be a continuous function, satisfying 
\[
f(\la) = 
\begin{cases}
\mathcal O(\la^{\tilde a + \tilde b+\tilde c - \tilde d-\hf-\eps}e^{-\pi \la}), & \quad \la \rightarrow \infty,\quad \eps>0,\\
\mathcal O(\la^\de), & \quad \la \downarrow 0, \quad \de>0.
\end{cases}
\]
Then
\[
\lim_{N \rightarrow \infty}\frac{1}{2\pi} \int_0^\infty f(\la) \inprod{\phi_\la}{ \phi_{\la'}}_N\, d\la=  \frac{f(\la')}{\tilde M(\la')}.
\]
\end{prop}
\begin{proof}
From the $c$-function expansion in Proposition \ref{prop:c-expa} we find
\[
[\phi_\la, \phi_{\la'}](N) = \sum_{\epsilon, \xi \in \{-1,1\}} 
\tilde c( \epsilon \la) \tilde c(\xi \la') \,[\Phi_{\epsilon\la}, \Phi_{\xi\la'}](N),
\]
and then we obtain from Lemma \ref{lem:asympt Wr} and Proposition \ref{prop:int Wr}
\[
\inprod{\phi_\la}{ \phi_{\la'}}_N = i\sum_{\epsilon, \xi \in \{-1,1\}} \tilde c(\epsilon \la) \tilde c( \xi \la')\,
\frac{N^{-2i( \epsilon \la+\xi \la')}}{ \epsilon \la + \xi \la'} 
\left(1 + \mathcal O \Big( \frac{1}{N} \Big) \right).
\]
We multiply both sides with an arbitrary function $f(\la)$, and we integrate over $\la$ from $0$ to $\infty$. The function $f$ must satisfy certain conditions that we determine later on. Letting $N \rightarrow \infty$ then gives
\[
\begin{split}
\lim_{N \rightarrow \infty} \int_0^\infty f(\la) \inprod{\phi_\la}{ \phi_{\la'}}_N&\, d\la=\\ 
\lim_{N \rightarrow \infty} i\int_0^\infty f(\la) \Big\{&
\psi_1(\la) \cos\big(2[\la+\la'] \ln(N) \big) +\psi_2(\la) \sin\big(2[\la+\la'] \ln(N) \big) \\ 
& + \psi_3(\la) \cos\big(2[\la - \la'] \ln(N) \big) + \psi_4(\la)
\frac{\sin\big(2[\la-\la']\ln(N)\big)}{\la-\la'} \Big\} d\la
\end{split}
\]
where
\[
\begin{split}
\psi_1(\la) &= \frac{1}{\la+\la'}\, \Big( \tilde c(\la) \tilde c(\la') - \tilde c(-\la) \tilde c(-\la') \Big), \\
\psi_2(\la) &= \frac{-i}{\la+\la'}\, \Big( \tilde c(\la) \tilde c(\la') + \tilde c(-\la) \tilde c(-\la') \Big), \\
\psi_3 (\la) & = \frac{1}{\la-\la'}\, \Big( \tilde c(\la) \tilde c(-\la') - \tilde c(-\la) \tilde c(\la') \Big), \\
\psi_4(\la) & = -i \Big( \tilde c(\la) \tilde c(-\la') + \tilde c(-\la) \tilde c(\la') \Big).
\end{split}
\]
Observe that $\tilde c(\la')  \tilde c(-\la') - \tilde c(-\la') \tilde c(\la')= 0$, so $\psi_3$ has a removable singularity at $\la = \la'$.
From the Riemann-Lebesgue lemma we find that the terms with $\psi_i$, $i=1,2,3$, vanish, provided that $f \psi_i \in L^1(0,\infty)$. We recognize the term with $\psi_4$ as a Dirichlet integral. Using the well-known property (see e.g. \cite[\S 9.7]{WW}) for Dirichlet integrals
\begin{equation} \label{eq:Dirichlet}
\lim_{t \rightarrow \infty}\frac{1}{\pi} \int_0^\infty g(x) \frac{
\sin[t(x-y)]}{x-y} dx = g(y),
\end{equation}
for a continuous function $g \in L^1(0,\infty)$, we obtain
\[
\lim_{N \rightarrow \infty} \int_0^\infty f(\la) \inprod{\phi_\la}{ \phi_{\la'}}_N\, d\la= \pi \psi_4(\la')f(\la')
=2\pi\,  \tilde c(\la')\tilde c(-\la') f(\la') = 2\pi \frac{f(\la')}{\tilde M(\la')}.
\]
The asymptotic behaviour of $\psi_i$, for $i=1,2,3,4$, can be obtained in the same way as the asymptotic behaviour of $\tilde M$, see \eqref{eq:asym M}. 
Then we find, for $i=1,2,3,4$,
\[
|\psi_i(\la)| = 
\begin{cases}
\mathcal O(\la^{-\hf-\tilde a-\tilde b -\tilde c + \tilde d}e^{\pi \la}),  &\la \rightarrow \infty, \\
\mathcal O(\la^{-1}), &\la \downarrow 0.
\end{cases}
\]
So, if $f$ satisfies the conditions given in the proposition, then $f\psi_i \in L^1(0,\infty)$. 
\end{proof}

Let $f$ be a continuous function. We define a linear operator $\mathcal F$ by
\[
(\mathcal F f)(\la) = \int f(x) \phi_\la(x) dm(x).
\]
We call $\mathcal F$ the Wilson function transform of type I. We denote the continuous part of the above integral by $\mathcal F_c f$, i.e.
\[
(\mathcal F_c f)(\la) = \frac{1}{2\pi}\int_0^\infty f(x) \phi_\la(x) M(x) dx.
\]

From the asymptotic behaviour of $\phi_\la$ and $M(x)$ we find that both $\mathcal Ff$ and $\mathcal F_c f$ are well defined if $f$ satisfies the conditions
\[
f(x) = 
\begin{cases}
\mathcal O(x^{d-a-b-c-\eps}e^{\pi x}), & x \rightarrow \infty, \quad \eps>0,\\
\mathcal O(x^{\de-1}), & x \downarrow 0, \quad \de>0.
\end{cases}
\]
Observe that if $f$ satisfies the condition given above, then $f \in \mathcal M$. Let $\mathcal M_0 \subset \mathcal M$ be the space consisting of continuous functions satisfying the above conditions. Then $\mathcal M_0$ is a dense subspace of $\mathcal M$ (it contains for instance the Wilson polynomials $R_n(x;a,b,c,1-d)$, which form an orthogonal basis for $\mathcal M$).

\begin{prop} \label{prop2}
Let $g \in \tilde{\mathcal M}_0$ and $\la \geq 0$, then
$\big(\mathcal F (\tilde{\mathcal F_c}g)\big)(\la) = g(\la)$.
\end{prop}
\begin{proof}
Define $f(\la) = \tilde M(\la) g(\la)$, then $f$ satisfies the conditions given in Proposition \ref{prop1}. Then we have
\[
\begin{split}
\big(\mathcal F (\tilde{\mathcal F_c}g)\big)(\la')=&
\int \phi_{\la'}(x)\Big( \frac{1}{2\pi}\int_0^\infty g(\la) \phi_\la(x) \tilde M(\la) d\la \Big) \ dm(x)\\
 =& \lim_{N\rightarrow \infty} \frac{1}{2\pi} \int_0^\infty f(\la) \Big(\frac{1}{2\pi} \int_0^N \phi_\la(x) \phi_{\la'}(x) M(x) dx +i\ \sum_{x \in \mathcal D}\ \phi_\la(x)\phi_{\la'}(x) \Res{z=x}\, M(z) \Big) \ d\la \\
=& \lim_{N\rightarrow \infty} \frac{1}{2\pi} \int_0^\infty f(\la) \langle \phi_\la, \phi_{\la'} \rangle_N d\la \\
=& \frac{f(\la')}{\tilde M(\la')} =  g(\la').
\end{split}
\]
Note that the first integral converges absolutely for $g \in \tilde{\mathcal M}_0$, so in that case interchanging the order of integration is allowed. 
\end{proof} 

In the next subsection we show that the dual Wilson function transform $\tilde{\mathcal F}$ is the inverse of the Wilson function transform $\mathcal F$. To do this, we must consider the discrete spectrum of the difference operator $L$.

\subsection{Discrete spectrum}
From the asymptotic behaviour of $\Phi_\la(x)$ and $M(x)$, see \eqref{eq:asympt} and \eqref{eq:asym M}, we obtain
\[
|\Phi_\la(x)|^2 M(x) = \mathcal O ( x^{4 \Im(\la)-1} ), \qquad x \rightarrow \infty.
\]
So for $\Im(\la)<0$ we have $\Phi_\la \in \mathcal M$. In this subsection we assume that $\la \in \tilde{\mathcal D}$ and that the set $\tilde{\mathcal D}$ is not empty, so $\Im(\la) <0$. In this case
$\tilde c(-\la)=0$, and therefore $\phi_\la(x) =\tilde c(\la) \Phi_\la(x) \in \mathcal M$. First we show that $\phi_\la$ is orthogonal to $\phi_{\la'}$ if $\la'\neq \la$.
\begin{prop} \label{prop:orth1}
For $\la \in \tilde{\mathcal D}$, $\la' \in \mathrm{supp}( d\tilde m)$, $\la \neq \la'$, we have
\[
\inprod{\phi_\la}{\phi_{\la'}}_{\mathcal M} = 0.
\]
\end{prop}
\begin{proof}
From Propositions \ref{prop:int Wr} and \ref{prop:c-expa} we obtain
\[
\inprod{\phi_\la}{\phi_{\la'}}_N = \frac{\tilde c(\la) \tilde c(\la') [\Phi_\la, \Phi_{\la'}](N)+ \tilde c(\la) \tilde c(-\la') [\Phi_\la, \Phi_{-\la'}](N)}{\la^2- \la'^2} .
\]
Then Lemma \ref{lem:asympt Wr} gives for large $N$ 
\begin{equation} \label{eq1}
\inprod{\phi_\la}{\phi_{\la'}}_N =   \frac{i(\la- \la') \tilde c(\la) \tilde c(\la')N^{-2i(\la+ \la')} +i(\la+ \la')\tilde c(\la) \tilde c(-\la')  N^{-2i(\la- \la')}}{ (\la+ \la')(\la-\la')} \left( 1+ \mathcal O\Big(\frac{1}{N} \Big) \right).
\end{equation}
Recall that for $\la \in \tilde{\mathcal D}$ we have $\la \in i\R_{<0}$. Then it is clear that for $\la'\in \R$ the right hand side tends to zero for $N \rightarrow \infty$. In case $\la' \in \tilde D$ the second term vanishes, and we have $\Im(\la + \la')<0$. So in this case the right hand side also tends to zero for $N \rightarrow \infty$.
\end{proof}
It remains to calculate the squared norm of $\phi_\la$ in case $\la \in \tilde {\mathcal D}$.
\begin{prop} \label{prop:orth2}
For $\la \in \tilde {\mathcal D}$
\[
\inprod{\phi_\la}{\phi_{\la} }_{\mathcal M}= \left(i\,\Res{\la'=\la}\ \tilde M(\la') \right)^{-1}.
\]
\end{prop}
\begin{proof}
We use expression \eqref{eq1}, where we let $\la' \rightarrow \la$. Then for large $N$
\[
\lim_{\la' \rightarrow \la}\inprod{\phi_\la}{\phi_{\la'}}_N =  \left(  \frac{i}{2\la } \tilde c(\la) \tilde c(\la) N^{4i \la} +i \tilde c(\la) \Big( \Res{\la'=\la}\, \frac{1}{\tilde c(-\la')} \Big)^{-1}  \right) \left( 1 + \mathcal O\Big(\frac{ 1}{N} \Big) \right).
\]
Letting $N \rightarrow \infty$ gives the result.
\end{proof}
\begin{rem}
If $\la, \la' \in \mathcal D_{\tilde a}$ then Propositions \ref{prop:orth1} and \ref{prop:orth2} give orthogonality relations for a finite number of Wilson polynomials with respect a measure that has only finitely many moments. Explicitly, if we assume $a, b, c, 1-d>0$, we have for $n,m <\hf(1-a-b-c-d)$
\[
\begin{split}
\frac{1}{2\pi}&\int_0^\infty R_n(x) R_m(x) \left| \frac{ \Ga(a+ix) \Ga(b+ix) \Ga(c+ix)}{ \Ga(1-d+ix) \Ga(2ix) } \right|^2 dx = \\ 
&\de_{nm}\frac{a+b+c+d-1}{a+b+c+d+2n-1} \frac{n! \, (b+c)_n (b+d)_n (c+d)_n} { (a+b)_n (a+c)_n (a+d)_n (a+b+c+d-1)_n } \\
&\times  \frac{\Ga(a+b) \Ga(a+c) \Ga(b+c) \Ga(1-a-b-c-d) }{  \Ga(1-a-d) \Ga(1-b-d) \Ga(1-c-d) },
\end{split}
\]
where $R_n$ is the Wilson polynomial defined by \eqref{def:Wilson pol}. Neretin \cite{Ne} found this orthogonality relation using an explicit evaluation of the Barnes-type integral
\[
\frac{1}{2\pi} \int_0^\infty \left| \frac{ \Ga(a+ix) \Ga(b+ix) \Ga(c+ix)}{ \Ga(1-d+ix) \Ga(2ix) } \right|^2 dx.
\]
Note that the above orthogonality relations remain valid for $n+m<1-a-b-c-d$.
\end{rem}

\subsection{The Wilson function transform: type I}
Combining Propositions \ref{prop:orth1} and \ref{prop:orth2} with Proposition \ref{prop2}, gives the following theorem.
\begin{thm} \label{Thm:WF I}
The Wilson function transform of type I, $\mathcal F : \mathcal M \rightarrow \tilde{\mathcal M}$,  defined by
\[
(\mathcal Ff)(\la) = \int f(x) \phi_\la(x) dm(x),
\]
is unitary, and its inverse is given by $\mathcal F^{-1} = \tilde{\mathcal F}$.
\end{thm}
\begin{proof}
First we show that $\mathcal F \circ \tilde{\mathcal F}$ is the identity operator on $\tilde{\mathcal M}_0$. For the continuous part of $\tilde{\mathcal F}$ this is Proposition \ref{prop2}, therefore we just write down the proof for the discrete part of $\tilde{\mathcal F}$. Let $g \in \tilde{\mathcal M}_0$. Recall that $\tilde{\mathcal D}$ is a finite set, then we obtain from Propositions \ref{prop:orth1} and \ref{prop:orth2}, 
\[
\begin{split}
\int \phi_{\la'}(x) \left( i\sum_{\la \in \tilde D} g(\la)\phi_\la(x) \Res{\mu=\la}\, \tilde M(\mu) \right) dm(x)
&= i\sum_{\la \in \tilde D} g(\la) \Res{\mu=\la}\, \tilde M(\mu) \left( \int \phi_\la(x) \phi_{\la'}(x)dm(x) \right)\\
& = 
\begin{cases}
0, & \la' \in [0,\infty),\\
g(\la') , & \la' \in \tilde D.
\end{cases}
\end{split}
\] 
Combining this with Proposition \ref{prop2} we obtain the desired result.
By duality we obtain $\big(\tilde{\mathcal F} (\mathcal Ff)\big)(x) = f(x)$ for $f \in \mathcal M_0$.

Next we show that $\tilde{\mathcal F}$ is an isometry on $\tilde{\mathcal M}_0$. For simplicity we assume that the measures $dm$ and $d\tilde m$ are absolutely continuous. From Proposition \ref{prop1} we obtain
\[
\begin{split}
\inprod{\tilde {\mathcal F} f}{\tilde {\mathcal F} g}_{\mathcal M} &=
\lim_{N \rightarrow \infty} \frac{1}{2\pi} \int_0^N (\tilde{\mathcal F} f)(x) \overline{(\tilde {\mathcal F} g)(x)}M(x) dx\\
& = \lim_{N \rightarrow \infty} \frac{1}{(2\pi)^2}\int_0^\infty \int_0^\infty f(\la) \overline{g(\la')} \inprod{\phi_\la}{\phi_{\la'}}_{N} \tilde M(\la) \tilde M(\la')\, d\la\, d\la'\\
&= \frac{1}{2\pi}\int_0^\infty f(\la') \overline{g(\la')} \tilde M(\la') d\la' = \inprod{f}{g}_{\tilde{\mathcal M}}.
\end{split}
\]
In case $dm$ and $d\tilde m$ have discrete mass points, the proof runs along the same lines. 

So the operator $\mathcal F:\mathcal M_0 \rightarrow \tilde{\mathcal M}_0$ is unitary. Since $\mathcal M_0$ is dense in $\mathcal M$, $\mathcal F$ extends uniquely to a unitary operator on $\mathcal M$.
\end{proof}

\begin{rem} An interesting special case is the case $a=b+c+d-1$. Then $(a,b,c,d) = (\tilde a, \tilde b, \tilde c, \tilde d)$, so the Wilson function transform of type I is then completely self-dual.
\end{rem}

\section{The Wilson function transform: type II} \label{sec II}
In this section we give another unitary integral transform which has the Wilson function as a kernel, by considering the action of the difference operator $L$ on (a dense subspace of) a different Hilbert space than in section \ref{sec:I}. The method we use is the same method as in section \ref{sec:I}, therefore we omit some details.

\subsection{The Hilbert space $\mathbf{\mathcal H}$}
Let $V^+ \subset \C^5$ be the set of parameters $a,b,c,d,t$ satisfying the following conditions:
\[
\overline{a}=1-d,\qquad  \overline{b} = c,\qquad  t \in \R.
\]
The dual parameters $\tilde a, \tilde b, \tilde c, \tilde d$ are still defined by \eqref{eq:dual param}, and we define the dual parameter $\tilde t$ by
\[
\tilde t = 1-\tilde b -c -t.
\]
It is easily verified that the assignment $(a,b,c,d,t) \mapsto (\tilde a, \tilde b, \tilde c, \tilde d, \tilde t)$ is an involution on $V^+$. Throughout this section we assume $(a,b,c,d,t) \in V^+$.

Let $H$ be the weight given by
\[
H(x)= H(x;a,b,c,d;t) = \frac{ \Ga(a \pm ix) \Ga(b \pm ix) \Ga(c \pm ix) \Ga(1-d \pm ix)}{ \sin\pi(t \pm ix)  \Ga(\pm 2ix)}.
\]
Here we use the notation $\sin(a \pm b) = \sin(a+b) \sin(a-b)$.
The weight $H$ is positive for $x \in \R$. Observe that $H(x;a,b,c,d;t) = M(x;a,b,c,d)/\sin\pi(t \pm ix)$, where $M$ is the weight function defined in section \ref{sec:I}. We assume that $H$ has only simple poles. This imposes conditions on the parameters that can be removed afterwards by continuity. Let $\mathcal D^+$ be the infinite discrete set defined by
\[
\mathcal D^+ = \{ i(t-n) \ | n \in \Z,\ t-n <0 \}.
\]
We define the measure $dh(\cdot)=dh(\cdot;a,b,c,d;t)$ by
\[
\int f(x) dh(x) = \frac{C}{2\pi}\int_0^\infty f(x) H(x) dx +iC \sum_{x \in \mathcal D^+}\ f(x) \Res{z=x}\, H(z).
\]
Here $C$ is the normalizing constant
\[
C = \sqrt{ \sin \pi(a+t) \sin \pi (b+t) \sin \pi(c+t) \sin \pi(1-d+t) }.
\]
Note that $\tilde C=C$.

For $x \in \mathcal D^+$ we have explicitly, 
\[
\begin{split}
i \Res{z=i(t-k)} \, H(z) = &  \Ga(a+t) \Ga(a-t) \Ga(b+t) \Ga(b-t) \Ga(c+t) \Ga(c-t) \Ga(1-d+t) \Ga(1-d-t)  \\
& \times \frac{k-t}{2t^2\pi^2} \frac{(a-t)_k (b-t)_k (c-t)_k (1-d-t)_k }{(1-a-t)_k (1-b-t)_k (1-c-t)_k (d-t)_k }
\end{split}
\] 
and then we see that for $(a,b,c,d,t) \in V^+$, the measure $dh$ is positive.
For $x\in \mathrm{supp}\ dh$ and $\la \in \mathrm{supp}\ d\tilde h$ the Wilson function $\phi_\la(x)$ is real valued.

We define the Hilbert space $\mathcal H=\mathcal H(a,b,c,d;t)$ to be the Hilbert space consisting of even functions that have finite norm with respect to inner product $\inprod{\cdot}{\cdot}_{\mathcal H}$ defined by
\[
\inprod{f}{g}_{\mathcal H} = \int f(x) \overline{g(x)} dh(x).
\]

\subsection{The Wronskian}
We denote $H_k= iC \Res{z = i(t-k)} H(z)$, and we define $k_0$ to be the smallest integer such that $t-k_0 < 0$. For $N \geq 0$ and $K \geq k_0$ we define a pairing $\langle \cdot, \cdot \rangle_{N,K}$ by
\[
\langle f, g \rangle_{N,K} =\frac{C}{2\pi}\int_0^N f(x)g(x) M(x) dx + \sum_{k=k_0}^K\ f\big(i(t-k)\big) g\big(i(t-k)\big) H_k.
\]
If $f$ and $g$ are real valued functions in $\mathcal H$, the limits $N,K \rightarrow \infty$ give the inner product $\langle f, g \rangle_{\mathcal H}$. We denote $\lim_{N \rightarrow \infty}\langle f, g \rangle_{N,K} = \langle f, g \rangle_{K}$, assuming that the limit exists. 

For analytic functions $f$ and $g$ we define the Wronskian $[f,g]$ by
\[
[f,g](k) = \Big\{ f\big(i(t-k)\big) g\big(i(t-k-1)\big) - f\big(i(t-k-1)\big) g\big(i(t-k)\big) \Big\} B\big(i(t-k)\big) H_k.
\]
\begin{prop}
For even analytic functions $f$ and $g$ we have
\[
\begin{split}
\langle L f, g \rangle_{N,K} - \langle  f,L g \rangle_{N,K} =&  \frac{C}{2\pi}\int_N^{N+i} \left\{ f(x) g(x-i) - f(x-i) g(x) \right\} B(x) H(x) dx -[f,g](K).
\end{split}
\]
\end{prop}
\begin{proof}
We follow the proof of Proposition \ref{prop: Wr}. For even functions $f$ and $g$ we have
\[
\langle f ,  g \rangle_{N,K} = \frac{1}{4\pi}\int_{\mathcal C_{N,K}} f(x) g(x) H(x) dx,
\]
where $\mathcal C_{N,K}$ is a contour in the complex plane defined as follows:
\begin{itemize}
\item $\mathcal C_{N,K}$ starts at $x=-N$ and ends at $x=N$,
\item $\mathcal C_{N,K}$ is invariant under reflection in the origin,
\item $\mathcal C_{N,K}$ separates the sequence of poles $i(a+n)$, $n \in \Z_{\geq 0}$, from the sequence $-i(a+m)$, $m \in \Z_{\geq 0}$, and similarly for poles of $H(x)$ corresponding to $b,c,1-d$, 
\item $C_{N,K}$ separates the sequence $i(t-n)$, $n \in\Z$, $n\leq K$, from the sequence $i(t-n)$, $n \in\Z$, $n \geq K+1$,
\item $\mathcal C_{N,K}$ separates the sequence of poles $i(t-n)$, $n \leq K$, from the sequence $-i(t-m)$, $m \leq K$, 
\end{itemize}
Using $B(-x+i)H(x-i) = B(x)H(x)$, we have
\[
\langle L f, g \rangle_{N,K} - \langle  f,L g \rangle_{N,K} =
\frac{C}{4\pi}\left(\int_{\mathcal C_{N,K}} -\int_{\mathcal C_{N,K}+i} \right) \Big( f(x-i) g(x) - f(x) g(x-i) \Big) B(x) H(x) dx.
\]
Now we make a counterclockwise oriented closed contour by connecting $\mathcal C_{N,K}$ and $\mathcal C_{N,K}+i$ at the end points, then the integrand $I(x)$ has two poles inside the closed contour; at $x= i(t-K)$ and $x=-i(t-K-1)$. So, 
\[
\frac{C}{4\pi}\left(\int_{\mathcal C_N} -\int_{\mathcal C_N+i}\right) I(x) dx =\frac{C}{4\pi } \left(\int_{-N}^{-N+i} - \int_N^{N+i}\right)I(x) dx + \frac{iC}{2} \Res{x=i(t-K)} I(x)+ \frac{iC}{2} \Res{x=-i(t-K-1)} I(x).
\]
In the same way as in the proof Lemma \ref{lem:oddWr} we can show that the first integral on the right hand side is equal to the second integral with opposite sign. For the residues we have
\[
\frac{iC}{2} \Res{x=i(t-K)} I(x) = -\hf[f,g](K),
\]
and, since $I(x) = -I(-x+i)$, 
\[
\begin{split}
\frac{iC}{2} \Res{x=-i(t-K-1)} I(x)& = -\lim_{x \rightarrow -i(t-K-1)} \big(x+i(t-K-1)\big)I(-x+i)\\
& = \lim_{y \rightarrow i(t-K)} \big(y-i(t-K)\big)I(y)= \Res{y=i(t-K)}I(y)=- \hf[f,g](K).
\end{split}
\]
This gives the desired result.
\end{proof}
From $H(x) = M(x) / \sin \pi(t \pm ix)$ and \eqref{eq:asym M} we find, for $y \in \R$,
\begin{equation} \label{eq:asym H}
H(x+iy) = \mathcal O(x^{2a+2b+2c-2d-1}e^{-4\pi x} ), \qquad x \rightarrow \infty.
\end{equation}
This gives the following for the Wilson functions.
\begin{prop} \label{prop:wron II}
For $\la \neq \la'$
\[
\langle \phi_\la, \phi_{\la'} \rangle_K = \frac{[\phi_\la, \phi_{\la'}](K) }{\la'^2- \la^2  }.
\]
\end{prop}
\begin{proof}
Since the Wilson function is an eigenfunction of $L$ for eigenvalue $\tilde a^2+ \la^2$, the result follows from Proposition \ref{prop:wron II}, because the integral on the right hand side in Proposition \ref{prop:wron II} is equal to zero. Indeed, using \eqref{eq:asym H} and the asymptotic behaviour of $\phi_\la$, which follows from Proposition \ref{prop:c-expa} and \eqref{eq:asympt}, and  applying dominated convergence, we obtain
\[
\lim_{N \rightarrow \infty} \int_N^{N+i} \left\{ \phi_\la(x) \phi_{\la'}(x-i) - \phi_\la(x-i) \phi_{\la'}(x) \right\} B(x) H(x) dx =0.
\]
\end{proof}
To find the the asymptotic behaviour of the Wronskian $[\phi_\la, \phi_{\la'}](K)$ for $K \rightarrow \infty$, we need the asymptotic behaviour of $H_K$, $B(i(t-K))$ and $\phi_\la(i(t-K))$. First we expand $\phi_\la$ in functions with nice asymptotic behaviour, as in Proposition \ref{prop:c-expa}.
\begin{prop} \label{prop:d-expan}
For $x = i(t-k) \in \mathcal D^+$
\[
\phi_\la(x) = \tilde d(\la) \Te_\la(k) + \tilde d(-\la) \Te_{-\la}(k),
\]
where
\[
\begin{split}
\Te_\la(k) &= \frac{(-1)^k}{\pi}\frac{ \Ga( 1-\tilde b-c-i\la -t+k)}{\Ga( \tilde b+c+i\la -t+k) } \F{4}{3}{\tilde a+i\la, \tilde b+i\la, \tilde c+i\la, 1-\tilde d+i\la}{ \tilde b+c+i\la-t+k, \tilde b+c+i\la +t-k, 1+2i\la}{1},\\
\tilde d(\la) &=  \frac{ \Ga(-2i\la) \sin\pi(\tilde t- i\la)}{ \Ga(\tilde a-i\la) \Ga(\tilde b-i\la) \Ga(\tilde c-i\la) \Ga(1-\tilde d-i\la)}.
\end{split}
\]
\end{prop}
\begin{proof}
This follows from Proposition \ref{prop:c-expa}, and Euler reflection formula;
\[
\frac{1}{ \Ga(\tilde b + c + i\la \pm ix) } = (-1)^k \frac{ \Ga(1-\tilde b -c -i\la + ix) \sin \pi(\tilde b+c+t+i\la)  }{ \pi\, \Ga(\tilde b +c +i\la +ix)}, .
\]
for $x = i(t-k)$, $k \in \Z$.
\end{proof}
\begin{rem}
Observe that 
\[
\frac{1}{\tilde d(\la) \tilde d(\la) } = \tilde H(\la).
\]
\end{rem}
For $k \rightarrow \infty$ we find from the explicit expression for $\Te_\la$ and \eqref{eq:asym Ga}, for $y \in \Z$,
\begin{equation} \label{eq:asym Te}
\Te_\la(k+y) = \frac{(-1)^k}{\pi} k^{d-a-b-c-2i\la} \left(1 + \frac{1}{2k} (d-a-b-c-2i\la)(1-2t+2y)+ \mathcal O \Big( \frac{1}{k^2}\Big) \right).
\end{equation}
Furthermore, for $k \rightarrow \infty$,
\begin{gather*}
H_k =  \frac{2\pi^2}{C} k^{2a+2b+2c-2d-1} \left(1+ \Big(\frac{1}{k}\Big)  \right),\\
B(i(t-k)) = \frac{k^2}{4}\left(1+ \Big(\frac{1}{k}\Big)  \right).
\end{gather*}

We can find the Wronskian $[\phi_\la, \phi_{\la'}](K)$, for $K \rightarrow \infty$, from the four Wronskians $[\Te_{\pm \la}, \Te_{\pm \la'}](K)$.
\begin{lem} \label{lem:asympt WrII}
For $K \rightarrow \infty$,
\[
[\Te_\la, \Te_{\la'}](K) = \frac{i}{C}(\la'-\la) K^{-2i(\la+ \la')}\left(1+ \Big(\frac1K \Big) \right).
\]
\end{lem}
\begin{proof}
We have
\[
[\Te_\la, \Te_{\la'}](K) = \Big\{\Te_\la(K) \Te_{\la'}(K+1) - \Te_\la(K+1) \Te_{\la'}(K)\Big\} B(i(t-K)) H_K.
\]
The lemma follows from this expression using the asymptotic behaviour of $\Te_\la(K)$, $B(i(t-K))$ and $H_K$.
\end{proof}

\subsection{Continuous spectrum}
In this subsection we assume $\la, \la' \geq 0$. 

\begin{prop} \label{prop1 II}
Let $f$ be an even continuous function, satisfying 
\[
f(\la) = 
\begin{cases}
\mathcal O(\la^{\tilde a + \tilde b+\tilde c - \tilde d-\hf-\eps}e^{-2\pi \la}), & \quad \la \rightarrow \infty,\quad \eps>0,\\
\mathcal O(\la^\de), & \quad \la \downarrow 0, \quad \de>0.
\end{cases}
\]
Then
\[
\lim_{N \rightarrow \infty}\frac{1}{2\pi} \int_0^\infty f(\la) \inprod{\phi_\la}{ \phi_{\la'}}_N\, d\la=  \frac{f(\la')}{C\,\tilde H(\la')}.
\]
\end{prop}
\begin{proof}
From Proposition \ref{prop:d-expan} we find,
\[
[\phi_\la, \phi_{\la'}](K) = \sum_{\epsilon, \xi \in \{-1,1\}} 
\tilde d( \epsilon \la) \tilde d(\xi \la') \,[\Phi_{\epsilon\la}, \Phi_{\xi\la'}](K),
\]
and then we obtain from Lemma \ref{lem:asympt WrII} and Proposition \ref{prop:wron II}
\[
\langle\phi_\la, \phi_{\la'}\rangle_K = \frac{i}{C}\sum_{\epsilon, \xi \in \{-1,1\}} \tilde d(\epsilon \la) \tilde d( \xi \la')\,
\frac{K^{-2i( \epsilon \la+\xi \la')}}{ \epsilon \la + \xi \la'} 
\left(1 + \mathcal O \Big( \frac{1}{K} \Big) \right).
\]
We multiply both sides with an arbitrary function $f(\la)$, and we integrate over $\la$ from $0$ to $\infty$. The function $f$ must satisfy certain conditions that we determine later on. Letting $K \rightarrow \infty$ then gives
\[
\begin{split}
\lim_{K \rightarrow \infty} \int_0^\infty f(\la) \inprod{\phi_\la}{ \phi_{\la'}}_K&\, d\la=\\ 
\lim_{K \rightarrow \infty} \frac{i}{C} \int_0^\infty f(\la) \Big\{&
\psi_1(\la) \cos\big(2[\la+\la'] \ln(K) \big) +\psi_2(\la) \sin\big(2[\la+\la'] \ln(K) \big) \\ 
& + \psi_3(\la) \cos\big(2[\la - \la'] \ln(K) \big) + \psi_4(\la)
\frac{\sin\big(2[\la-\la']\ln(K)\big)}{\la-\la'} \Big\} d\la
\end{split}
\]
where
\[
\begin{split}
\psi_1(\la) &= \frac{1}{\la+\la'}\, \Big( \tilde d(\la) \tilde d(\la') - \tilde d(-\la) \tilde d(-\la') \Big), \\
\psi_2(\la) &= \frac{-i}{\la+\la'}\, \Big( \tilde d(\la) \tilde d(\la') + \tilde d(-\la) \tilde d(-\la') \Big), \\
\psi_3 (\la) & = \frac{1}{\la-\la'}\, \Big( \tilde d(\la) \tilde d(-\la') - \tilde d(-\la) \tilde d(\la') \Big), \\
\psi_4(\la) & = -i \Big( \tilde d(\la) \tilde d(-\la') + \tilde d(-\la) \tilde d(\la') \Big).
\end{split}
\]
$\psi_3$ has a removable singularity at $\la = \la'$. From the Riemann-Lebesgue lemma we find that the terms with $\psi_i$, $i=1,2,3$, vanish, provided that $f \psi_i \in L^1(0,\infty)$. Using \eqref{eq:Dirichlet} for the term with $\psi_4$, we obtain
\[
\begin{split}
\lim_{N \rightarrow \infty} \int_0^\infty f(\la) \inprod{\phi_\la}{ \phi_{\la'}}_N\, d\la= \frac\pi C \psi_4(\la')f(\la')
=\frac{2\pi}{C}  \tilde d(\la')\tilde d(-\la') f(\la') = \frac{ 2\pi}{C} \frac{f(\la')}{\tilde H(\la')}.
\end{split}
\]
Using, for $i=1,2,3,4$,
\[
|\psi_i(\la)| = 
\begin{cases}
\mathcal O(\la^{-\hf-\tilde a-\tilde b -\tilde c + \tilde d}e^{2\pi \la}),  &\la \rightarrow \infty, \\
\mathcal O(\la^{-1}), &\la \downarrow 0,
\end{cases}
\]
we see that if $f$ satisfies the conditions given in the proposition, then $f\psi_i \in L^1(0,\infty)$. 
\end{proof}

Let $\mathcal H_0$ be the dense subspace of $\mathcal H$ defined by
\[
\mathcal H_0 = \left\{ f \in \mathcal H \, | \, f\big(i(t-k)\big)=0\ \text{for} \ k \gg 0 \right\}.
\]
For $f \in \mathcal H_0$ we define a linear operator $\mathcal G$ by
\[
(\mathcal G f)(\la) = \int f(x) \phi_\la(x) dh(x).
\]
We call $\mathcal G$ the Wilson function transform of type II. We denote the continuous part of the above integral by $\mathcal G_c f$, i.e.
\[
(\mathcal G_c f)(\la) = \frac{C}{2\pi}\int_0^\infty f(x) \phi_\la(x) H(x) dx.
\]
From the asymptotic behaviour of $\phi_\la$ and $H(x)$ we find that both $\mathcal Gf$ and $\mathcal G_c f$ are well defined if $f\in \mathcal H_0$. 
\begin{prop} \label{prop2 II}
Let $g$ be a continuous function satisfying 
\[
g(\la) = 
\begin{cases}
\mathcal O\big(\la^{\tilde d- \tilde a - \tilde b - \tilde c + \hf - \eps} e^{2\pi \la} \big), & \la \rightarrow \infty, \quad \eps>0,\\
\mathcal O\big(\la^{\de-1}), & \la \downarrow 0, \quad \de>0,
\end{cases}
\] 
then $\big(\mathcal F (\tilde{\mathcal F_c}g)\big)(\la) = g(\la)$.
\end{prop}
\begin{proof}
The proof runs along the same lines as the proof of Proposition \ref{prop2}.
\end{proof} 
Note that if $g \in \tilde{\mathcal H}_0$, then $g=o(\la^{\tilde d - \tilde a-\tilde b-\tilde c} e^{2\pi \la})$ for $\la \rightarrow \infty$. So if $g \in \mathcal H_0$, then $g$ satisfies the conditions of Proposition \ref{prop2 II}. 

\subsection{Discrete spectrum} In this subsection we assume $\la \in \tilde{\mathcal D}^+$. In this case $\tilde d(-\la) =0$ and $\Im\la<0$, and then it follows from Proposition \ref{prop:d-expan} and the asymptotic behaviour \eqref{eq:asym Te} of $\Te_\la$, that $\phi_\la \in \mathcal H$. The following two propositions are proved in a similar way as Propositions \ref{prop:orth1} and \ref{prop:orth2}, therefore we omit the proofs here.
\begin{prop} \label{prop:orth1 II}
For $\la \in \tilde{\mathcal D}^+$, $\la' \in \mathrm{supp}( d\tilde h)$, $\la \neq \la'$, we have
\[
\inprod{\phi_\la}{\phi_{\la'}}_{\mathcal H} = 0.
\]
\end{prop}
\begin{prop} \label{prop:orth2 II}
For $\la \in \tilde{\mathcal D}^+$,
\[
\langle \phi_\la, \phi_{\la}\rangle_{\mathcal H} =  \left( iC\, \Res{\la'=\la}\, \tilde H(\la') \right)^{-1}.
\]
\end{prop}

\subsection{The Wilson function transform: type II}
In the same way as in the proof of Theorem \ref{Thm:WF I} we can 
combine Propositions \ref{prop:orth1 II} and \ref{prop:orth2 II} with Proposition \ref{prop2 II}, to find that the operator $\mathcal G: \mathcal H_0 \rightarrow \tilde{\mathcal H}_0$ is a unitary operator with inverse $\tilde{\mathcal G}$. Since $\mathcal H_0$ is dense in $\mathcal H$, $\mathcal G$ extends uniquely to a unitary operator on $\mathcal H$.
\begin{thm}
The Wilson function transform of type II, $\mathcal G : \mathcal H \rightarrow \tilde{\mathcal H}$, defined by
\[
(\mathcal G f)(\la) = \int f(x) \phi_\la(x) dh(x),
\]
is unitary, and its inverse is given by $\mathcal G^{-1} = \tilde{\mathcal G}$.
\end{thm}
\begin{rem}
For $a=b+c+d-1$, $t=\hf(1-b-c-s)$, $s\in \R$, the Wilson function transform of type II is completely self-dual.
\end{rem}

\section{Explicit transformations} \label{sec:expl}
In this section we calculate explicitly the Wilson function transforms of certain functions. First we give an integral representation for the Wilson function related to Jacobi functions. This leads to two explicit transformations in Theorems \ref{thm:trJ} and \ref{thm:tr4F3}. Then, in Theorem \ref{thm:trans Wpol}, we show that the Wilson function transform of type I maps an orthogonal basis of polynomials in the Hilbert space $\mathcal M$ to itself, with dual parameters.

\subsection{Transformations related to Jacobi functions}
The Jacobi functions, see \cite{Koo}, are defined by
\begin{equation} \label{def:Jac funct}
\varphi_\la^{(\al,\be)}(x) = \F{2}{1}{\hf(\al+\be+1-i\la), \hf(\al+\be+1+i\la)}{\al+1}{-x}.
\end{equation}
For $x \geq 1$ we use here the unique one valued analytic continuation of the $_2F_1$-function. The Jacobi functions are the kernel in an integral transform pair, called the Jacobi transform, given by
\begin{equation} \label{def: Jacobi trans}
\begin{cases}
\displaystyle(\mathcal J f)(\la)= \int_0^\infty f(z)\varphi_\la^{(\al,\be)}(x) \De_{\al,\be}(x) dx \\ \\
\displaystyle f(x) = \frac{1}{2\pi}\int (\mathcal J f)(\la)\varphi_\la^{(\al,\be)}(x) d\nu(\la)
\end{cases}
\end{equation}
where $\al>-1$, $\be \in \R \cup i\R$,
\[
\De_{\al,\be}(x)= x^{\al}(1+x)^{\be},
\]
and $d\nu(\la)$ is the measure given by
\[
\begin{split}
\frac{1}{2\pi}\int g(\la) d\nu(\la) &= \frac{1}{2\pi}\int_0^\infty g(\la) |c_{\al,\be}(\la)|^{-2} d\la -i \sum_{\la \in \mathcal E} g(\la) \Res{\mu=\la} \big(c_{\al,\be}(\mu) c_{\al,\be}(-\mu) \big)^{-1}, \\
c_{\al,\be}(\la)&= \frac{2^{-i\la} \Ga(\al+1) \Ga(i\la) }{ \Ga\big(\hf(\al+\be+1+i\la)\big) \Ga\big(\hf(\al-\be+1+i\la)\big) },\\
\mathcal E &= \Big\{ i(|\be|-\al-1-2j)\ | \ j \in \Z_{\geq 0}, |\be|-\al-1-2j>0 \Big\}.
\end{split}
\]

The Jacobi polynomials $P_n^{(\al,\be)}(x)$ are orthogonal polynomials on the interval $[-1,1]$ with respect to the measure $(1-x)^\al (1+x)^\be dx$. They  have an explicit expression as $_2F_1$-series, see \cite{AAR}, \cite{KS}. In this paper we need Jacobi polynomials of argument $ (1-x)/(1+x)$.  Using Pfaff's transformation, we find that the explicit expression for these polynomials is 
\begin{equation} \label{def:Jac}
P_n^{(\al,\be)}\left(\frac{1-x}{1+x}\right) = \frac{ ( \al+1)_n }{n!} (1+x)^{-n} \F{2}{1}{-n, -\be -n }{2\al}{-x}.
\end{equation}
The orthogonality relation for these polynomials reads, for $\Re(\al), \Re(\be)>-1$,
\begin{equation} \label{eq:orthJac}
\begin{split}
\int_0^\infty P_n^{(\al,\be)}\left( \frac{1-x}{1+x} \right)  & P_n^{(\al,\be)}\left( \frac{1-x}{1+x} \right) x^\al (1+x)^{-\al-\be -2} dx = \\ &\frac{1}{2n+ \al + \be+1 } \frac{ \Ga(n+\al+1) \Ga(n+\be+1) }{n!\, \Ga(n+ \al + \be +1)}\, \de_{nm}.
\end{split}
\end{equation}

The following proposition gives a representation of a very-well poised $_7F_6$-series as an integral over a product of two $_2F_1$-series.
\begin{prop} \label{prop:Wil-Jac}
For $\be,\mu\in \C$, $\Re(\al)>0$, $\Re(\ga+\de \pm \rho)>0$, 
\begin{equation} \label{eq:Wil-Jac}
\begin{split}
\int_0^\infty & \F{2}{1}{\al+\mu+\ga,\al+\mu-\ga}{ 2\al}{-u} 
\F{2}{1}{\al+\be+\rho,\al+\be-\rho}{ 2\al}{-u} u^{2\al-1} (1+u)^{\be-\de+\mu} du\\
=& \frac{\Ga(2\al) \Ga(2\al + \de +\rho+\ga) \Ga(\de+\ga\pm \rho) \Ga(\de-\ga+\rho)} {\Ga(\al+\de+\rho \pm \mu) \Ga(\al+\de+\ga \pm \be)}\\
&\times W(2\al+\de+\ga-1+\rho; \al+\ga+\mu, \al+\ga-\mu, \al+\be+\rho, \al-\be+\rho,  \de+\ga+\rho).
\end{split}
\end{equation}
\end{prop}
\begin{proof}
We start with the following formula, for $\Re(\al)>0$, $\Re(\ga+\de \pm \rho)>0$, \begin{equation} \label{eq:h1}
\begin{split}
\int_0^1 y^{2\al-1} (1-y)^{\ga+\de-\al-\be-1} &\F{2}{1}{\al+\be+\rho, \al+\be-\rho}{2\al}{\frac{y}{1-y}}dy\\
 =& \frac{ \Ga(2\al) \Ga(\ga+\de + \rho) \Ga(\ga+\de - \rho)}{ \Ga(\al-\be+\ga+\de)\Ga(\al+\be+\ga+\de)}.
\end{split}
\end{equation}
To prove this identity we transform the $_2F_1$-series by Pfaff's transformation \cite[Thm.2.2.5]{AAR}, then we obtain a $_2F_1$-series that converges uniformly on $[0,1]$ for $\Re(\rho)<0$. We interchange summation and integration, then the result follows from using the Beta-integral and Gauss's summation formula \cite[Thm.2.2.2]{AAR}. The condition on $\rho$ can be removed using the symmetry in $\rho$ and $-\rho$, and continuity in $\rho$.

Next we write the integral in the theorem as
\[
\begin{split}
 I = \int_0^1 & y^{2\al-1} (1-y)^{\de-\mu-\be -2\al-1} \\ \times& \F{2}{1}{\al+\mu+\ga,\al+\mu-\ga}{ 2\al}{\frac{y}{y-1}} 
\F{2}{1}{\al+\be+\rho,\al+\be-\rho}{ 2\al}{\frac{y}{y-1}}  dy,
\end{split}
\]
where we substituted $u \mapsto y/(1-y)$. By \cite[\S2.10(3)]{Erd} the first $_2F_1$-function can be expanded in terms of $_2F_1$-functions of argument $1-y$
\begin{equation} \label{eq:2F1+2F1}
\begin{split}
&\F{2}{1}{\al+\mu+\ga, \al+\mu-\ga }{2\al}{\frac{y}{y-1}} = \\
&\frac{ \Ga(2\al) \Ga(-2\ga)}{\Ga(\al+\mu-\ga) \Ga( \al-\mu-\ga)} (1-y)^{\al+\mu+\ga } \F{2}{1}{\al+\ga +\mu , \al+\ga-\mu }{1+2\ga}{ 1-y }\\
 +& \frac{ \Ga(2\al) \Ga(2\ga)}{\Ga(\al+\mu+\ga) \Ga( \al-\mu+\ga)} (1-y)^{\al+\mu-\ga } \F{2}{1}{\al-\ga +\mu , \al-\ga-\mu }{1-2\ga}{ 1-y }.
\end{split}
\end{equation}
Observe that the second term is equal to the first term with $\ga$ replace by $-\ga$. The integral $I$ splits according to this as $I=I_\ga+I_{-\ga}$.  We use formula \eqref{eq:h1} to evaluate $I_\ga$.
Interchanging summation and integration, which is allowed for $\Re(\al)<\hf$ since then the $_2F_1$-series converges uniformly on $[0,1]$, and using \eqref{eq:h1} leads to 
\begin{equation} \label{eq:4F3}
\begin{split}
I_\ga=&\frac{ \Ga(2\al)^2 \Ga(-2\ga) \Ga(\de+\ga +\rho) \Ga(\de+\ga -\rho) }{\Ga(\al-\ga+ \mu) \Ga(\al-\ga- \mu) \Ga(\al+\be+\ga+\de) \Ga(\al-\be+\ga+\de)} \\
&\times \F{4}{3}{ \al+\ga+\mu, \al+\ga-\mu, \de+\ga+\rho, \de+\ga-\rho }{2\ga+1, \al+\be+\ga+\de, \al-\be+\ga+\de}{1}.
\end{split}
\end{equation}
Then $I=I_\ga+I_{-\ga}$ is the sum of two balanced $_4F_3$-functions, which by \cite[\S4.4(4)]{Bai} can be written as the very-well poised $_7F_6$-series given in the proposition. 

Note that the expression in \eqref{eq:4F3} is an analytic function in $\al$ for $\Re(\al) >0$. The integrand of the integral $I$ is analytic in $\al$ for $\Re(\al)>0$, and continuous in $y$. So differentation with respect to $\al$ and integration with respect to $y$ can be interchanged, see e.g. \cite[\S4.2]{WW}. We see that the integral $I$ is an analytic function in $\al$ for $\Re(\al)>0$, and therefore the condition $\Re(\al)<\hf$ can be removed by analytic continuation.
\end{proof}

Both $_2F_1$-series in Proposition \ref{prop:Wil-Jac} can be considered as Jacobi functions. Using the substitutions
\begin{equation} \label{eq:subst}
\al \mapsto \hf(a+1-d), \quad \be \mapsto \hf(c-b), \quad \ga \mapsto \hf(a+d-1), \quad \de \mapsto \hf(b+c), \quad \mu \mapsto ix, \quad \rho \mapsto i\la,
\end{equation}
and the definition \eqref{def:Wilson function} of the Wilson functions, the right hand side of \eqref{eq:Wil-Jac} can be written as 
\[
\Ga(1+a-d)^2 \Ga(\tilde a \pm i\la) \Ga(1-\tilde d \pm i\la)\, \phi_\la(x;a,b,c,d).
\]
So Proposition \ref{prop:Wil-Jac} gives a representation of the Wilson function as  the Jacobi function transform of a Jacobi function. From the inverse Jacobi transform we obtain the following.
\begin{thm} \label{thm:trJ} 
Let $\varphi(x) = \varphi_{2x}^{(b+c-1,c-b+1)}(u)$, $u\geq 0$, and $\psi(x) = \sin \pi(t \pm ix) \varphi(x)$, then
\[
\begin{split}
(\mathcal F \varphi)(\la) &=
(1+u)^{\hf-\tilde c+i\la} \F{2}{1}{\tilde a+i\la, 1-\tilde d+i\la}{ 1+\tilde a - \tilde d }{-u},\qquad (a,b,c,d) \in V,\\
(\mathcal G \psi)(\la) &= C\,(1+u)^{\hf-\tilde c+i\la} \F{2}{1}{\tilde a+i\la, 1-\tilde d+i\la}{ 1+\tilde a - \tilde d }{-u},\qquad (a,b,c,d,t) \in V^+.
\end{split}
\]
\end{thm}
Observe that for $u=0$ the first statement gives $(\mathcal F 1)(\la)=1$.
\begin{proof}
In Proposition \ref{prop:Wil-Jac} we replace $\rho$ by $i\la \in i\R$. Then using the definition of the Jacobi functions \eqref{def:Jac funct} we can write Proposition \ref{prop:Wil-Jac} as
\[
\int_0^\infty (1+u)^{\mu-\de-\be}\F{2}{1}{\al+\mu+\ga,\al+\mu-\ga}{ 2\al}{-u} \varphi_{2\la}^{(2\al-1, 2\be)}(u) \De_{2\al-1, 2\be}(u) du = F(\la),
\]
where $F(\la)$ denotes the right hand side of \eqref{eq:Wil-Jac}. Note that from \eqref{eq:2F1+2F1} we find for $u \rightarrow \infty$
\[
(1+u)^{\mu-\de-\be}\F{2}{1}{\al+\mu+\ga,\al+\mu-\ga}{ 2\al}{-u} \sim C_1 u^{-\al-\be-\de - \ga} + C_2u^{-\al-\be-\de+ \ga},
\]
where $C_1$ and $C_2$ are independent of $u$. So for $\Re(\de \pm \ga)>0$ we see that this function is an element of $L^2([0,\infty), \De_{2\al-1,2\be}(u)du)$. 
Taking the inverse Jacobi transform, assuming for simplicity that the discrete set $\mathcal E$ is empty, gives
\[
\begin{split}
(1+u)^{\mu-\de-\be}&\F{2}{1}{\al+\mu+\ga,\al+\mu-\ga}{ 2\al}{-u} = \\
&\frac{1}{2\pi} \int_0^\infty F(\la) \varphi_{2\la}^{(2\al-1, 2\be)}(u) \left| \frac{ \Ga(\al+\be + i\la) \Ga( \al-\be+i\la) }{ \Ga(2\al) \Ga(2i\la)} \right|^2 d\la.
\end{split}
\]
We write out $F(\la)$ and $\varphi_{2\la}^{(2\al-1,2\be)}(t)$ as hypergeometric functions and use the substitutions given in \eqref{eq:subst}, without the last substitution, then we obtain
\[
\begin{split}
(1+u)^{\hf-c+ix} \F{2}{1}{a+ix,1-d+ix}{ 1+a-d}{-u} =
 \frac{1}{2\pi} 
\int_0^\infty  \F{2}{1}{\tilde c +i\la, \tilde c -i\la }{\tilde b +\tilde c}{-u} \phi_\la(x) \tilde M(\la) d\la .
\end{split}
\]
This is the first statement in the theorem with dual parameters, and $x$ and $\la$ interchanged. Writing the above integral as a contour integral (with the contour as in the proof of Proposition \ref{prop: Wr} with $N \rightarrow \infty$), and using analytic continuation in the parameters $a,b,c,d$, the identity can be extended to the case $(a,b,c,d) \in V$. Deforming the contour again to the real line might add a finite number of discrete mass points.

The second statement follows from the first. We use $H(x) = M(x) / \sin\pi(t \pm ix)$ and we use that the function $\sin \pi(t \pm ix) \varphi(x)$ vanishes on the infinite discrete set $\mathcal D^+$.
\end{proof}

We substitute $\mu = - \al - \ga  -n $, $n \in \Z_{\geq 0}$, in Proposition \ref{prop:Wil-Jac}, then the first $_2F_1$-series on the left hand side of \eqref{eq:Wil-Jac} terminates and can be written as a Jacobi polynomial by \eqref{def:Jac}. Writing the right hand side as $I_\ga+I_{-\ga}$, with $I_\ga$ as in \eqref{eq:4F3}, we see that $I_{-\ga}=0$ because of the factor $\Ga(\al+\ga+\mu)^{-1}$. Now we obtain, for $\Re(\al)>0$ and $\Re(\ga+\de \pm \rho)>0$,
\begin{equation} \label{eq:Koo}
\begin{split}
\int_0^\infty u^{2\al-1} (1+u)^{\be-\de-\al-\ga} & P_n^{(2\al-1,2\ga)}\left( \frac{ 1-u }{1+u} \right) \F{2}{1}{\al+\be+\rho,\al+\be-\rho}{ 2\al}{-u} du \\
=&\frac{ (-1)^n (2\ga+1)_n \Ga(2\al) \Ga(\ga+\de + \rho) \Ga(\ga+\de - \rho) }{ n!\, \Ga(\al + \be +\ga + \de )  \Ga(\al - \be +\ga + \de )} \\
\times &\F{4}{3}{ -n, 2\al+2\ga+n, \de+\ga+\rho, \de+\ga-\rho }{2\ga+1, \al+\be+\ga+\de, \al-\be+\ga+\de}{1}.
\end{split}
\end{equation}
This is Koornwinder's formula \cite[(3.3)]{Koo1} stating that Jacobi polynomials are mapped onto Wilson polynomials by the Jacobi function transform. We use this formula to expand the integral in Proposition \ref{prop:Wil-Jac} in terms of terminating $_4F_3$-functions. This leads to an expansion formula for the Wilson function.

\begin{prop} \label{prop:expan}
Let $f,g \in \C$ satisfy $\Re(f)>0$, $\Re(f+g)>0$, $\Re(\tilde a +g)>0$, $\Re(1-\tilde d +g)>0$. Then 
\[
\begin{split}
\phi_\la(x;a,b,c,d) = \sum_{n=0}^\infty &C_n(x,\la) \F{4}{3}{-n, a+f+g-d+n, f+i\la, f-i\la}{f+g, \tilde b +f, \tilde c+f}{1} \\
& \times \F{4}{3}{ -n, a-d+f+g+n, g+ \tilde a, g +1-\tilde d }{ f+g, g+ \tilde b +c +ix, g+ \tilde b +c -ix}{1},
\end{split}
\]
where
\[
\begin{split}
C_n(x,\la) =& \frac{ \Ga (a-d+f+g+1) \Ga(g+ \tilde a) \Ga(g+1-\tilde d) \Ga( f \pm i\la)}{ \Ga(1+a-d) \Ga(f+g) \Ga(f+ \tilde b) \Ga(f+\tilde c)\Ga(g + \tilde b +c \pm ix) \Ga( \tilde a \pm i\la) \Ga(1-\tilde d \pm i\la)}\\
&\times \, \frac{ a+f+g-d+2n}{ a+f+g-d}\, \frac{ (a+f+g-d)_n (f+g)_n }{ n!\, (a+1-d)_n}.
\end{split}
\]

\end{prop}
\begin{proof}
Let $I$ be the integral in Proposition \ref{prop:Wil-Jac}. We expand the first $_2F_1$-function in $I$ in terms of Jacobi polynomials of argument $(1-u)/(1+u)$:
\[
(1+u)^{\al+\mu+ \eta - \si +1} \F{2}{1}{ \al+\mu+ \ga, \al+\mu-\ga}{2\al}{-u} = \sum_{n=0}^\infty c_n \, P_n^{(2\al-1,2\eta)} \left( \frac{ 1-u }{1+u} \right) .
\]
Note that from \eqref{eq:2F1+2F1} we obtain that the left hand behaves for large $t$ as $C_1 u^{\eta -\si - \ga+1} + C_2 u^{\eta -\si + \ga+1}$, where $C_1$ and $C_2$ are independent of $t$. So for $\Re(\eta+\si \pm \ga)>0$ the function on the left hand side is an element of $L^2\big( (0,\infty), u^{2\al-1}(1+u)^{-2\al-2\eta-1}du\big)$. We see that the expansion on the right hand side converges uniformly for $u$ in compact intervals. 
The coefficients $c_n$ are found using \eqref{eq:Koo} and the orthogonality relation for the Jacobi polynomials \eqref{eq:orthJac}, and this gives for $\Re(\al)>0$, $\Re(\eta)>-\hf$, and $\Re(\si+\eta \pm \ga)>0$,
\[
\begin{split}
c_n =& \frac{\Ga(2\al+2\eta+1) \Ga(\eta+\si \pm \ga)  }{\Ga(2\eta+1) \Ga(\al  +\eta + \si \pm \mu) }\ \frac{(2\al+2\eta+2n)}{(2\al+2\eta)} \frac{ (-1)^n (2\al+2\eta)_n}{ (2\al)_n}  \\
&\times \F{4}{3}{ -n, 2\al+2\eta+n, \si+\eta+\ga, \si+\eta-\ga }{2\eta+1, \al+\mu+\eta+\si, \al-\mu+\eta+\si}{1}
\end{split}
\]
From this expansion we find 
\[
I= \sum_{n=0}^\infty c_n \int_0^\infty u^{2\al-1} (1+u)^{\be-\al-\de-\eta+\si-1} P_n^{(2\al-1,2\eta)} \left( \frac{ 1-u }{1+u} \right) \F{2}{1}{\al+\be+\rho,\al+\be-\rho}{ 2\al}{-u} du
\]
Using \eqref{eq:Koo} again gives the following expansion for $I$,
for $\Re(\al)>0$, $\Re(\eta)>-\hf$, $\Re(\si+\eta \pm \ga)>0$, $\Re(\de-\si+\eta+1 \pm \rho)>0$, 
\[
\begin{split}
I=&\frac{\Ga(2\al) \Ga(2\al+2\eta+1) \Ga(\eta+\si \pm \ga) \Ga(\eta+\de-\si+1 \pm \rho)  }{\Ga(2\eta+1) \Ga(\al  +\eta + \si \pm \mu) \Ga(\al  +\eta + \de-\si+1 \pm \be ) }\\
\times &\sum_{n=0}^\infty 
\frac{2\al+2\eta+2n}{2\al+2\eta}\, \frac{  (2\al+2\eta)_n (2\eta+1)_n}{  n! \, (2\al)_n}  \F{4}{3}{ -n, 2\al+2\eta+n, \si+\eta+\ga, \si+\eta-\ga }{2\eta+1, \al+\mu+\eta+\si, \al-\mu+\eta+\si}{1}\\
&\qquad \times \F{4}{3}{ -n, 2\al+2\eta+n, \de-\si+1+\eta+\rho, \de-\si+1+\eta-\rho }{2\eta+1, \al+\be+\eta+\de-\si+1, \al-\be+\eta+\de-\si+1}{1}.
\end{split}
\]
The proposition then follows from Proposition \ref{prop:Wil-Jac}, the substitutions $\eta-\si + 1 + \de \mapsto f$, $\eta + \si-\de \mapsto g$ and \eqref{eq:subst}.
\end{proof}

\begin{cor} \label{cor:expan}
For $\Re(f)>0$
\[
\begin{split}
\phi_\la&(x;a,b,c,d) =\\ 
&\sum_{n=0}^\infty C_n(x,\la) \F{3}{2}{-n, f+i\la, f-i\la}{ \tilde b +f, \tilde c +f}{1} \F{3}{2}{ -n , f-\tilde a, f+\tilde d -1}{ f + \tilde c -c -ix, f + \tilde c -c +ix }{1},
\end{split}
\]
where
\[
C_n(x,\la) = \frac{(f+\tilde c -c +ix)_n (f+\tilde c -c -ix)_n \, \Ga(f \pm i\la) }{n!\,\Ga(1+a-d+n) \Ga(f+\tilde b) \Ga(f+\tilde c) \Ga(\tilde a \pm i\la) \Ga(1-\tilde d \pm i\la)}.
\]
\end{cor}
\begin{proof}
We use Whipple's transformation \cite[Thm.3.3.3]{AAR} to write the second $_4F_3$-function in Proposition \ref{prop:expan} as
\[
\frac{ (f+\tilde c -c +ix)_n (f+\tilde c -c -ix)_n } {(g+\tilde b +c +ix)_n (g+\tilde b +c +ix)_n } \F{4}{3}{ -n, a-d+f+g+n, f-\tilde a, f+\tilde d -1 }{ f+g, f+\tilde c -c -ix, f+ \tilde c -c +ix}{1}.
\]
Then the corollary follows from letting $g \rightarrow \infty$, using \eqref{eq:asym Ga} for the $\Ga$-functions.
\end{proof}

We consider one of the $_4F_3$-functions in Proposition \ref{prop:expan} as a Wilson polynomial, then we obtain the following Wilson function transform of a Wilson polynomial.
\begin{thm} \label{thm:tr4F3}
Let $f,g\in \R$ such that $g+a$, $g+b$, $g+c$, $g+1-d$ has positive real part. Furthermore,  let $\varphi_n$ and $\psi_n$ be the function defined by
\[
\varphi_n(x) = \Ga(g \pm ix) R_n(x;f,b,c,g), \quad \psi_n(x)= \sin \pi(t \pm ix)\varphi_n(x).
\]
Then 
\[
\begin{split}
(\mathcal F \varphi_n )(\la)= C^{-1}(\mathcal G \psi_n)(\la)&=\frac{ \Ga(g + a ) \Ga(g + b ) \Ga(g+c) \Ga(g+1- d)}{ \Ga(g+\tilde b +c \pm i\la) } \frac{ (g + b)_n (g+ c)_n }{ (f+ b)_n (f+ c)_n}\\
&\times  \F{4}{3}{ -n , b+c+f+g+n-1, g+1- d, g+ a}{ f+g, g+\tilde b +c+i\la, g+\tilde b +c -i\la}{1}.
\end{split}
\]
\end{thm}
Here we assume for $\mathcal F$ that $(a,b,c,d) \in V$, and for $\mathcal G$ that $(a,b,c,d,t) \in V^+$.

\begin{proof}
We recognize the first $_4F_3$-function in Proposition \ref{prop:expan} as the Wilson polynomial $R_n(\la)=R_n(\la;f,\tilde b, \tilde c, g)$. We multiply by $R_n(\la)$ and integrate against the orthogonality measure $d\mu(\la;f,\tilde b, \tilde c, g)$. Using the orthogonality relation \eqref{eq:orth} then gives
\[
\begin{split}
\frac{1}{2\pi} \int_0^\infty &R_n(\la) \phi_\la(x;a,b,c,d) \frac{ \Ga(\tilde a \pm i\la) \Ga( \tilde b \pm i\la) \Ga( \tilde c \pm i\la) \Ga(1-\tilde d \pm i\la) \Ga(g \pm i\la) }{ \Ga( \pm 2i\la) }\, d\la= \\
&\frac{ \Ga(g +\tilde a ) \Ga(g +\tilde b ) \Ga(g+\tilde c) \Ga(g+1-\tilde d) }{ \Ga(g+\tilde b +c \pm ix) } \frac{ (g + \tilde b)_n (g+ \tilde c)_n }{ (f+\tilde b)_n (f+\tilde c)_n}\\
& \times \F{4}{3}{ -n , a-d+f+g+n, g+1-\tilde d, g+\tilde a}{ f+g, g+\tilde b +c +ix, g+\tilde b+c -ix}{1}.
\end{split}
\]
The theorem follows from this by replacing the parameters $a,b,c,d$ by their dual parameters, and interchanging $\la$ and $x$. 

Note that the integral can again be written as a contour integral. Then by analytic continuation the result remains true when $g \in \R$ and the pairwise sum of $a,b,c,d,g$ has positive real part.
\end{proof}
From Corollary \ref{cor:expan} a similar theorem can be derived for the Wilson function transform of a continuous dual Hahn polynomial, which is a limit case of the Wilson polynomial, see \cite{AAR}, \cite{KS}. We leave this to the reader.

\subsection{The Wilson function transform of type I of an orthogonal system}
Koornwinder's formula \eqref{eq:Koo} shows that the Jacobi function transform maps a (bi)orthogonal system of polynomials (the Jacobi polynomials) in $L^2\big((0,\infty), \De_{\al, \be}(t)dt\big)$, onto an orthogonal system of polynomials (the Wilson polynomials) in $L^2\big((0,\infty), d\nu\big)$. We show that there exists a similar formula for the Wilson function transform of type I.\\ 

The Wilson polynomials $R_n(x) = R_n(x;a,b,c,1-d)$ form an orthogonal basis for the Hilbert space $\mathcal M$. We show that the Wilson function transform of type I maps a Wilson polynomial $R_n$  onto itself, with dual parameters. The proof is based on the following proposition. 
\begin{prop} \label{prop:LR}
Let $L$ be the difference operator as in Proposition \ref{prop:L}, then
\[
(L R_n)(x) = C_n R_{n+1}(x) + (\tilde C_n + \tilde D_n)R_n(x)+ D_n R_{n-1}(x),
\]
where 
\[
\begin{split}
C_n &= \frac{ (n+a+b+c-d)(n+a+b)(n+a+c)(n+a+1-d) }{(2n+a+b+c-d)(2n+a+b+c-d+1)}, \\
D_n &= \frac{ n (n+b+c-1)(n+b-d)(n+c-d) }{(2n+a+b+c-d)(2n+a+b+c-d-1)}.
\end{split}
\]
\end{prop}
\begin{proof}
Applying the difference operator $L$ on $R_n$ gives
\[
\begin{split}
(L R_n )(x)& = B(-x) R_n(x+i) + \big[A(-x)+A(x)\big]R_n(x) + B(x) R_n(x-i) \\
&= \big[ n(n+a+b+c-d) + A(x)+A(-x)+B(x) + B(-x) \big] R_n(x).
\end{split}
\]
Here we use the difference equation \eqref{eq:diff} for the Wilson polynomials $R_n(x;a,b,c,1-d)$. From the explicit expressions for $A(x)$ and $B(x)$ we obtain
\[
A(x)+B(x) = \frac{(a+ix)(b+ix)(c+ix) }{2ix},
\]
and this gives
\[
\begin{split}
A(x)+A(-x)+B(x)+B(-x)& = \frac{(a+ix)(b+ix)(c+ix)-(a-ix)(b-ix)(c-ix) }{2ix}\\
&= ab+ac+bc-x^2.
\end{split}
\]
Then, using the three-term recurrence relation \eqref{eq:rec Wpol} for the Wilson polynomials, we find
\[
\begin{split}
(L R_n)(x) =& \big[ n(n+a+b+c-d) + ab +ac +bc +a^2 \big]R_n(x) - (a^2+x^2)R_n(x) \\
=& C_n R_{n+1}(x)  +D_n R_{n-1}(x)\\
&+[ n(n+a+b+c-d) + ab +ac +bc +a^2  -C_n -D_n] R_n(x).
\end{split}
\]
Note that we use here $C_n$ and $D_n$ as in \eqref{eq:rec Wpol}, but with $d$ replaced by $1-d$. A long, but straightforward calculation shows that 
\[
C_n +D_n + \tilde C_n + \tilde D_n = n(n+a+b+c-d) + a^2+ab+ac+bc,
\]
and from this the proposition follows.
\end{proof}

\begin{thm} \label{thm:trans Wpol}
For $(a,b,c,d) \in V$, the Wilson function transform of type I of the Wilson polynomial $R_n(x)=R_n(x;a,b,c,1-d)$ is given by
\[
(\mathcal F R_n)(\la) = (-1)^n \frac{ (b+c)_n }{(1+a-d)_n} \tilde R_n(\la).
\]
\end{thm}
\begin{proof}
We use $(\mathcal F R_n)(\la)= \lim_{N \rightarrow \infty} \langle \phi_\la,R_n \rangle_{N }$, where $\langle \cdot, \cdot \rangle_N$ is the pairing defined in section \ref{sec:I}. From Proposition \ref{prop: Wr} we find, 
\[
\lim_{N \rightarrow \infty} \langle L \phi_\la,R_n \rangle_{N } - \langle \phi_\la,L R_n\rangle_{N } = \lim_{N \rightarrow \infty} [\phi_\la, R_n](N).
\]
Using \eqref{eq:asym M}, \eqref{eq:asympt}, $R_n(x)= \mathcal O(x^{2n})$, and applying dominated convergence, we obtain for the Wronskian 
\[
\lim_{N \rightarrow \infty}[\phi_\la, R_n](N) =0.
\]
 The Wilson function $\phi_\la$ is an eigenfunction of $L$ for eigenvalue $\tilde a^2 + \la^2$, so we have
\[
(\tilde a^2+\la^2)(\mathcal F R_n)(\la) = \lim_{N \rightarrow \infty} \langle L \phi_\la,R_n\rangle _{N} =\lim_{N \rightarrow \infty} \langle  \phi_\la,LR_n\rangle _{N}= \big(\mathcal F (LR_n)\big)(\la).
\]
Then from Proposition \ref{prop:LR} and linearity of $\mathcal F$, it follows that the function $\mathcal F R_n$ satisfies the recurrence relation
\begin{equation} \label{eq:rec rel}
(\tilde a^2+\la^2)\,y_n(\la) = C_n\, y_{n+1}(\la) + (\tilde C_n + \tilde D_n )\, y_n(\la) +D_n\, y_{n-1}(\la).
\end{equation}
From $R_{-1}(x)=0$, we obtain $(\mathcal F R_{-1})(\la)=0$. By Theorem \ref{thm:trJ} we have $(\mathcal F R_0)(\la)=(\mathcal F 1)(\la)=1$, so we find from the recurrence relation that $(\mathcal F R_n)(\la) = P_n(\la)$. Here $P_n(\la)$ is a polynomial in $\la^2$ of degree $n$ satisfying the three-term recurrence relation \eqref{eq:rec rel}, with initial values $P_{-1}(\la)=0$ and $P_0(\la)=1$. From the recurrence relation for $\tilde R_n(\la)= R_n(\la;\tilde a, \tilde b, \tilde c, 1-\tilde d)$ it follows that $P_n(\la) = (-1)^n \frac{ (b+c)_n }{(1+a-d)_n} \tilde R_n(\la)$. 
\end{proof}
\begin{rem}
Theorem \ref{thm:trans Wpol} gives in fact a new proof for the unitarity of $\mathcal F$, since an orthogonal basis in $\mathcal M$ is mapped to an orthogonal basis in $\tilde{\mathcal M}$ with the same norm. Note that Theorem \ref{Thm:WF I} is not used to proof Theorem \ref{thm:trans Wpol}.
\end{rem}

Theorem \ref{thm:trans Wpol} can be considered as the 
$q=1$ analogue of \cite[Prop.4.1]{St}, which is used by Stokman to obtain an expansion of the Askey-Wilson function in Askey-Wilson polynomials, see \cite[Thm.4.2]{St}. So one might expect that from Theorem \ref{thm:trans Wpol} a similar expansion can be proved for the Wilson functions. However, if we formally expand
\[
\phi_\la(x;a,b,c,d) = \sum_{n=0}^\infty d_n(x,\la)\, R_n(x;a,b,c,1-d) R_n(\la;\tilde a, \tilde b, \tilde c, 1-\tilde d),
\]
and we calculate $d_n(x,\la)$ using Theorem \ref{thm:trans Wpol} and the orthogonality relation for the Wilson polynomials, we find an expansion that 
in general does not converge, but it does converge if $x \in \mathcal D$, or $\la \in \tilde{\mathcal D}$. The convergence of the expansion formula for the Askey-Wilson functions is due to a Gaussian factor $q^{m(m+1)/2}$, $0<q<1$. In the limit $q \rightarrow 1$ this factor disappears.

From Theorem \ref{thm:trans Wpol} we can also find the Wilson function transform of type II of a Wilson polynomial.
\begin{thm}
For $(a,b,c,d,t) \in V^+$, the Wilson function transform of type II of the Wilson polynomial $R_n(x)=R_n(x;a,b,c,1-d)$ is given by
\[
(\mathcal G R_n)(\la) = (-1)^n \frac{ (b+c)_n }{(1+a-d)_n} \frac{\sin \pi(\tilde t \pm i\la) }{ \sqrt{\sin \pi(a+t) \sin \pi(b+t) \sin \pi(c+t) \sin \pi(1-d+t) }} \ \tilde R_n(\la).
\]
\end{thm}
\begin{proof}
From Theorem \ref{thm:trans Wpol} we obtain
\[
\big(\mathcal G \, \sin \pi(t \pm ix) R_n\big)(\la) =  (-1)^n  \frac{C\, (b+c)_n }{(1+a-d)_n} \tilde R_n(\la).
\]
Taking the inverse of this gives
\[
(\tilde{\mathcal G} \tilde R_n)(x) = (-1)^n \frac{ (1+a-d)_n }{ C \, (b+c)_n} \sin \pi (t \pm ix) R_n(x).
\]
Replacing all parameters by their dual and interchanging $x$ and $\la$ then gives the statement in the theorem. 
\end{proof}

\end{document}